\documentclass[10pt]{article}
\usepackage[sumlimits]{amsmath}
\usepackage{amsmath,amsfonts,amsthm,amssymb,amscd}
\usepackage[latin1]{inputenc}  
\usepackage[dvipsnames]{xcolor}

\newtheorem{theo}{Theorem}[section]
\newtheorem{prop}[theo]{Proposition}
\newtheorem{lemm}[theo]{Lemma}
\newtheorem{coro}[theo]{Corollary}

\newtheorem{defi}[theo]{Definition}

\def\bbb{{\cal B}}\def\ddd{{\cal D}}
\def\fff{{\cal F}}
 \def\hhh{{\cal H}}   

\def\lll{{\cal L}}
\def\nnn{{\cal N}} 
  
\def\sss{{\cal S}}\def\ttt{{\cal T}}

\newcommand\N{\mathbb{N}}

\newcommand\R{\mathbb{R}}
\newcommand\C{\mathbb{C}}

\def\1{{\bold 1}}
\def\ph{\varphi}
\def\eps{{\varepsilon}}

\title{Positivity results for  Weyl's pseudodifferential calculus on
the Wiener space.}
\author{L. Jager }

\begin{document}

\maketitle

\begin{abstract}
This paper deals with positivity properties for  a
pseudodifferential calculus, generalizing Weyl's
classical quantization, and set 
on an infinite dimensional phase space, the Wiener space. \\
In this frame, we show that a positive symbol does not, in general, give
a positive operator.
In order to measure the nonpositivity, we  establish a G\aa rding's
inequality, which holds for the symbol classes at hand. 
Nevertheless, for symbols with radial aspects, additional assumptions ensure
the positivity of the associated operator. 
  
\end{abstract}


\section{Introduction and main results}\label{sec-intro}

This paper comes after a series of articles defining  a pseudodifferential
calculus on the Wiener space, which is an infinite dimensional measure space.
This calculus, constructed and developed
in 
\cite{AJN-JFA, AJN-NMR,A-L-N-2,J-StExt}, 
generalizes Weyl's well-known calculus.
The aim of the construction was to treat problems from mathematical physics,
in which an unpredictable number of particles may appear. This explains the
necessity of replacing the finite dimensional phase and configuration spaces
by infinite dimensional Hilbert spaces.
An argument in favour of  pseudodifferential calculus is that it allows
working with functions (the symbols) rather than dealing with  operators.
This  makes Weyl's calculus, sometimes,  preferable  to the use of the
Fock space
(which is linked with the Wiener space anyway).\\

The preceding articles contain different constructions of the
calculus, $L^2$-boundedness properties, the definition of two different
symbol classes.
Tools and results analogous to the finite dimensional
results exist, such as a Beals characterization
\cite{A-L-N-2,AJN-NMR}, composition results \cite{AJN-comp},
in the shape of semiclassical asymptotic expansions in
powers of a small parameter $h$. One of the constructions relies on 
the notion of Wigner function, as in the finite dimensional case.
Parallel calculus are available, like the Anti-Wick calculus
(associating an operator with a symbol) and the Wick 
calculus (associating a function with an operator). A special kind of
heat operators links these calculi together \cite{J-StExt, AJN-comp,AJN-NMR}.
The Wick calculus
is defined thanks to a family of coherent states, which are
elements of a space of square summable functions defiened on the Wiener space
and have a counterpart in the Fock space.\\

In this paper, we add to this theory by proving positivity and nonpositivity
results.
 Some answer natural questions : if the symbol is nonnegative,
  what can we say about the operator ? Is Flandrin's conjecture (about integrals of Wigner functions on convex sets) valid ?
Proposition \ref{posit-rad-gene} establishes that, under particular conditions, a
positive symbol gives a positive operator. Under
more general conditions, a positive symbol gives
an  operator which is not too negative (G\aa rding's inequality, Proposition \ref{Gaa}). \\

The paper is organized as follows.
Before stating the main results, the 
introduction gives the necessary definitions about Weyl's pseudodifferential
calculus on the Wiener space. In particular, we 
recall the definition of the Wiener space, which consists of a 
Hilbert space $\hhh$ and of a  normed probability space $B$ containing $\hhh$
and replaces the configuration space $\R^n$. 
This is followed by the constructions of the calculus, one using Wigner
functions, the other one  relying on symbol classes. A
$L^2$-boundedness result is recalled after.
To keep this first section as light as possible, other notions are  delayed
until  Section \ref{sec-Gaarding}, in which they are ingredients of the proofs.\\

Section \ref{sec-Ext-stoch} contains some complements about stochastic
extensions in particular cases, mainly for the  so-called cylindrical
functions. Stochastic extensions are the main means of turning a function
defined on $\hhh$ into a function defined on $B$, when it is possible. It is
a probabilistic definition, the topological usual extensions or restrictions
failing to  operate. This section gives the existence and an explicit
expression for stochastic extensions in the case, for example, of Wigner
functions. \\

Section \ref{sec-pos} gives the first results about positivity. They are
presented together since the arguments are rather similar. They intensely use 
 the Wigner functions, either finite dimensional or not, 
and a decomposition on a Hilbert basis (of $L^2(B)$) consisting of
Hermite functions. Section \ref{sec-Ext-stoch}  is
mainly  exploited here, since stochastic extensions of  cylindrical functions
allow the transposition of finite dimensional results in an infinite
dimensional frame. For example, Weyl's calculus on the Wiener space
is no more positive than its finite dimensional
predecessor and model. The  Flandrin conjecture does not
hold either (as was proved in finite dimension in \cite{DDL}). We give a
result about positivity in the case of radial symbols or of tensor products of
radial symbols too. \\

Finally, section \ref{sec-Gaarding} proves G\aa rding's inequality for
this calculus. 
It is more or less independent of Sections \ref{sec-Ext-stoch} and
\ref{sec-pos} and relies
on notions and results from \cite{AJN-JFA} which are recalled here,
like partial heat operators and hybrid Weyl-Anti-Wick operators. Although
shorter than section \ref{sec-pos}, this part may probably  find more
applications in the sequel. \\

{\it Conventions and notations.} \\
In this work we  denote by  $\bbb(X)$ the Borel $\sigma$-algebra  of
a normed space $X$ and by $\fff(X)$ the set of all  finite dimensional
subspaces of $X$. On a (finite dimensional) Euclidean space $E$, we
denote by $\mu_{E,s}$ the Gaussian measure of variance $s>0$. If $E=\R^n$,
$\mu_{\R^n,s}= (2\pi s)^{-n/2} e^{-\frac{|x|^2}{2s}} dx$, where $dx$ is the
Lebesgue measure.
\\
We have tried to denote differently the functions defined on the Hilbert
spaces $\hhh,\hhh^2$ and the functions defined on the Wiener extensions
$B$ or $B^2$.  The latter wear a $\tilde{\ }$ or are indexed by $B$.
Similarly, notions
relative to two different spaces, for example relative
either to $\hhh$ or to $B$,  are indexed by $\hhh$ or $\bbb$. This is
the case for  Wigner functions. \\

The remaining paragraphs of this section now recall the notions which
allow us  to state  the
results.
\subsection{The Wiener space}\label{sec-Wiener-space}
Most of the notions recalled here are defined in 
\cite{K} (chap. 1 par. 4), which sums up the construction of the articles
\cite{G-1, G-2, G-3, G-4}. Stochastic extensions appear in \cite{RA} too.\\

Let $\hhh$ be a real, separable, infinite dimensional Hilbert space,
with norm $| \ | $ or $| \ |_{\hhh}$ associated with the scalar
product $\cdot$. It would be practical to endow $\hhh$ with a measure
similar to the Lebesgue measure or a Gaussian measure, for example.
This is not directly possible  and one needs to enlarge $\hhh$ in
a convenient way.\\

The first step is to define a protomeasure on particular subsets of $\hhh$
called {\it cylinders}. Let $E \in \fff(\hhh)$ be a finite dimensional
subspace of $\hhh$, let $P_E$ be the orthogonal projection of $\hhh$ on $E$.
A cylinder based on $E$ is the inverse image, by $P_E$, of a Borel set $A$
of $E$. Explicitly, the cylinder defined by 
$$ C= \{ x\in \hhh \ : \ P_E(x)\in A \}$$
has a protomeasure given by 
\begin{equation}\label{proto-measure}
\mu_{\hhh,s}(C)= \int_A e^{-\frac{|y|^2}{2s}  }
(2\pi s)^{-\dim(E)/2} dy
= \int_A d\mu_{E,s}(y).
\end{equation}
We have implicitly chosen an orthonormal basis on $E$, $dy$ is the Lebesgue
measure on $E$ and the positive parameter $s$
represents the variance. The number  $\mu_{\hhh,s}(C)$   does not depend on
the space $E$ on which the cylinder is based (nor on the choice of the basis).
Similarly, a {\it cylindrical } or {\it tame} function on $\hhh$ is a function
$f$
which can be written as  $f=\ph\circ P_E$ for a given $E\in \fff(\hhh)$
and a function $\ph$ defined on $E$.  One may think of cylindrical functions
as depending on a finite
number of variables. \\

The protomeasure  \eqref{proto-measure} is finitely additive on the set of
cylinders but it is not $\sigma$-additive (unless one restricts oneself
to cylinders based on  a fixed $E$).  Therefore, it can't be extended as a
measure
on the $ \sigma$-algebra generated by the cylinders. To obtain a measure,
one introduces a new and larger set containing $\hhh$. Let $|| \ ||$
be another norm  on $\hhh$, satisfying the  condition below, classically called
{\it measurability}:
\begin{equation}\label{measurability}\forall\eps>0,\ \exists E_{\eps}\in \fff(\hhh) :\
\forall F\in \fff(\hhh),\ F\bot E_{\varepsilon} ,\
\mu_{\hhh,s}(\{x\in \hhh : ||P_{F}(x)||>\varepsilon\}) < \varepsilon.
\end{equation}
For this new norm, all $d$-dimensional subspaces $F$ of $\hhh$ are not on
the same level. With the original norm, the cylinder
$(\{x\in \hhh : |P_{F}(x)|>\varepsilon\})$ has a protomeasure  $\mu_{\hhh,s}$
which depends only on ${\rm dim}(F)$ and $\eps$, but not on the situation of
$F$ with respect to a space $E_{\eps}$. \\

One denotes by $B$ the completion of $\hhh$ with respect to $||\ ||$. It
is called a {\it Wiener extension } of $\hhh$, it depends on the choice of
the norm $|| \ ||$. The dual space of $B$
is called $B'$ and $\hhh$ is identified with its dual space. This gives
the sequence of inclusions:
$$B'\subset \hhh\subset B,$$
where each space is a dense subset of the following one and the inclusions are
continuous. The couple  $(\hhh,B)$ (norm and inclusions remaining implicit) is
called an {\it abstract Wiener space}. The new norm on $\hhh$
is not necessarily
associated with a scalar product and even if it were, it would not be
equivalent to the first one. \\

One may then define a measure on the cylinders of the Wiener extension $B$. For 
 $y_1,\dots, y_n\in B'$ and $A$ a Borel set of $ \R^n$, set:  
\begin{equation}\label{real-measure}
\mu_{B,s}(\{ x\in B \ : \ ((y_i,x)_{B',B})_{1\leq i\leq n} \in A\})
=\mu_{\hhh,s}(\{ x\in \hhh \ : \ (y_i\cdot x)_{1\leq i\leq n} \in A\}).
\end{equation}

This expression is, formally, analogous to the definition of the protomeasure,
but it gives a real probability measure on the $\sigma$-algebra
generated by the cylinders of $B$. This  $\sigma$-algebra
is the Borel $\sigma$-algebra of $B$,
which is a consequence of the separability of $\hhh$ and of $B$.
Note that $\hhh$ is dense in $B$ but negligible for $\mu_{B,s}$.
Even if we do not use the following fact in the paper, we may mention that the
symmetrized Fock space $\fff_s(\hhh)$ is isometric to any 
$L^2(B, \mu_{B,s})$.\\

We now define fundamental simple functions which will play an important part
in this paper. They replace the first degree monomials in the finite
dimensional theory and appear, for example, in Section \ref{basedim1}, as
elementary components of a Hilbert basis of a space $L^2(B)$.
Right now they allow one to define the notion of stochastic extension
defined, for example, in \cite{RA}.\\

An element $a$ of $B'$ can be seen in three different ways.
Of course, it is a linear
continuous form on $B$. Since $B'\subset \hhh$, it gives a linear
form from $\hhh$ to $\R$, identified with the form 
$x \mapsto x\cdot a$ defined on $\hhh$.
The point is that, since $B$ is endowed with the
$\sigma$-algebra $\bbb(B)$ and the probability measure
$\mu_{B,s}$, this element $a$ is also a {\it random variable} on
$(B,\bbb(B),\mu_{B,s})$. To stress the difference of status,
we denote the random
variable by $\ell_a$. Definition \eqref{real-measure} of the measure
implies that $\ell_a$ has the Gaussian distribution
$\nnn(0, \sigma^2=s|a|_{\hhh}^2)$ and
 its norm in $L^2(B,\mu_{B,s})$ is equal to $\sqrt{s}|a|_{\hhh}$.
By density (of $B'$ in $\hhh$) we obtain an application
from  $\hhh$ into $L^2(B,\mu_{B,s})$: 
\begin{equation}\label{application-ell}
\ell : \hhh\ \rightarrow \ L^2(B,\mu_{B,s}),\quad a\  \mapsto\  \ell_a. 
\end{equation}
Remark that
$\frac{1}{\sqrt{s}} \ell $ is isometric. The set of all $\ell_a, a\in \hhh$ is a
Gaussian Hilbert space in the sense of \cite{Jan}.\\

If $a$ belongs to $\hhh \setminus B'$, $\ell_a$ is only a random variable
and there is no
reason why it should be linear. Strictly speaking, it is defined only
almost everywhere. Still, it satisfies the relationships
$\ell_{a}(x+y)= \ell_a(x) + a\cdot  y$ for
$y\in \hhh$ and  $\ell_{a}(-x)= -\ell_a(x) $.
Sometimes we will write  $\ell_{a+ib}$ instead of
$\ell_a + i\ell_b$, for  $a,b\in \hhh$
(recall that  $\hhh$ is a real space). \\

One generalizes the projections, replacing the scalar product with an element
$a$ by the corresponding function $\ell_a$. Precisely, if 
$(u_j)_{j\leq \dim(E)} $ is a Hilbert basis of $E\in \fff(\hhh)$,
the orthogonal projection $P_E$ (from $\hhh$ to $E$)  and its extension
 $\tilde{\pi}_E$ (from $B$ to  $E$) are written below: 
\begin{equation}\label{projection}
  \forall x\in \hhh, \
  P_E (x)=\sum_{j=1}^{\dim(E)} (x\cdot u_j)  u_j,\quad 
  \forall y\in B, \
  \tilde{\pi}_E (y)=\sum_{j=1}^{\dim(E)} \ell_{u_j}(y)  u_j.
\end{equation}
The operators  $\tilde{\pi}_E$ are defined on $B$, whereas the
orthogonal projections are defined on $\hhh$. If
$E\in\fff(B')$,  $\tilde{\pi}_E$ is defined  everywhere on $B$ (and not
just ``almost everywhere'') and linear, since  
$ \ell_{u_j}=u_j\in B' $ for all $j$. \\

The generalized projections allow one to extend, in a certain sense,
functions initially defined on $\hhh$:

\begin{defi}\label{ext-stoch} 
  Let $( \hhh,B)$  be an abstract Wiener space.
Let $s$ be a positive real number.
\\
A function $f$ defined on $\hhh$  is said to admit a stochastic extension 
 $\widetilde f\in L^p (B,\mu _{B,s})$ in the sense of $L^p (B,\mu _{B,s})$ 
($1\leq p< \infty$) if, for every increasing sequence $(E_n)_{n\in \N}$ in
${\cal F}(\hhh)$,
 whose union is dense in $\hhh$, the functions  $f \circ\widetilde \pi_{E_n}$
 are in  $L^p (B,\mu_{B ,s})$ and if the sequence  $f \circ \widetilde \pi_{E_n}$
 satisfies
 $$f \circ \widetilde \pi_{E_n} \quad \longrightarrow \quad \widetilde f
\quad {\text in }\quad L^p (B,\mu _{B,s}).$$
\end{defi}

The original notion of Gross and Ramer required a convergence in probability,
which is implied by the $L^p$ convergences above.
In general, a stochastic  extension is not a continuity extension, though it is sometimes the case  (Theorem 6.3 \cite{K}).
In the same way, restricting to $\hhh$ a function defined on $B$ does
not always make sense, since $\hhh$ is negligible for all measures 
$\mu_{B,s}$, $s >0$.\\

A function defined on $\hhh$ has not necessarily a stochastic extension, or
it may have an extension which is useless. For example, 
$x\in \hhh \mapsto |x|^2$ or $x\mapsto e^{i|x|^2}$ have no
stochastic extension. The first example justifies the introduction of the
notion of a measurable norm \eqref{measurability}, see
\cite{K}, Chap.1, Sec. 4.
The function
$f : x\mapsto e^{-|x|^2}$ admits, as an extension, the null function on $B$.


\subsection {Weyl's calculus in the Wiener space} \label{sec-Weyl-Wiener}
In this section we recall the main definitions pertaining to the
infinite dimensional Weyl calculus.
Since there are two spaces $\hhh$ and $B$ and two different constructions,
we felt that a short guideline stressing the main features could be useful. 
\\

The {\it first  construction}, which is the less technical one, associates,
with a symbol $\tilde{F}$  defined on $B^2$, a quadratic or bilinear form
$Q_h^W(\tilde{F})$, which applies to a couple of convenient cylindrical
functions (see Definition \ref{space-D} for the meaning of convenient,
Definition
\ref{Weyl-quadratic-form} for the form $Q_h^W(\tilde{F})$).
This form is defined by an integral on $B^2$ and corresponds to the
following expression
in the finite dimensional situation : 
\begin{equation}\label{Weyl-classique-Wigner}
<Op_h^{Weyl,cl}(F)u,v>_{L^2(\R^n, d\lambda(x))}=  (2\pi h)^{-n}  
\quad \displaystyle
\int_{\R^{2n}} F(z,\zeta)  H^{cl}_h(u,v)(z,\zeta) dz d\zeta.
\end{equation}
Here,  $ H^{cl}_h(u,v)$ is the finite dimensional Wigner function of the couple
$(u,v)$. \\

The {\it second construction} introduces  {\it symbol classes}, as in the
finite
dimensional Weyl calculus.  The symbols belonging to these classes
are defined on $\hhh^2$ and satisfy regularity conditions, their partial
derivatives or Fréchet differential are bounded in a precise way 
(Definitions \ref{CV-Class}, \ref{Classe-DD-def}).
In this frame, one associates, with a symbol $F$ defined on $\hhh^2$,
an operator  $Op_h^{W}(F)$, which is linear and bounded on a space
$L^2(B)$. The operator is not defined directly, it is the limit of a Cauchy
sequence of rather complicated hybrid operators, which are recalled later on in
Section \ref{sec-Gaarding} but are not necessary right now to state the result.
The convergence takes place in the space of
the operators
bounded on a $L^2(B)$.
Contrary to the first construction,
 $Op_h^{W}(F)u$ has no integral expression generalizing
\begin{equation}\label{Weyl-classique}
 (Op_h^{Weyl,cl }(F) (u ))(x)  = (2 \pi h)^{-n} \int
_{\R^{2n}} e^{\frac{i}{h} (x - y)\cdot \xi } F \left ( \frac{x+y}{2} , \xi
\right )u (y) dy d\xi .
\end{equation}
When the operator is applied to (convenient) cylindrical functions, there
is a link with
the form 
$Q_h^W(\tilde{F})$, where $\tilde{F}$ is a stochastic extension of $F$. \\

These notions are taken from 
\cite{AJN-JFA,J-StExt,JG-rev} and are now  recalled at  length.\\

\noindent{\it First construction }\\

Let $E$ be a $d$-dimensional Euclidean space, identified with $\R^d$ by the
choice of an orthonormal basis. For $h>0$,
one defines an isometric isomorphism between 
$L^2(E,\mu_{\R^d,h/2})$ and $L^2(E,dy )$, setting
\begin{equation}\label{gamma}
\forall y\in E,\quad \gamma_{E,h/2} f(y)=
(\pi h)^{-d/4} e^{-\frac{|y|^2}{2h}}f(y).
\end{equation}
The space of the ``test'' functions, to which the quadratic form will be
applied, is given by the following definition.
\begin{defi}\label{space-D}
Let  $E$ belong to $\fff(B')$.
\begin{itemize}
 \item
One denotes by $\sss_{E,h/2}$  the space of all functions
$\ph$ defined on $E$, such that $\gamma_{E,h/2} \ph$ is rapidly decreasing.
\item
One denotes by $\ddd_{E,h/2}$ the set of all functions $\tilde{f}$ defined on
$B$ and
based on $E$, of the form   $\tilde{f} = \ph \circ \widetilde{\pi}_E$,
with $\ph\in \sss_{E,h/2}$ and 
$\widetilde{\pi}_E$ as in  \eqref{projection}.
\item
One then sets
$$
\ddd_{B',h/2}= \bigcup_{E\in \fff(B')} \ddd_{E,h/2}, \quad
\ddd_{\hhh,h/2}= \bigcup_{E\in \fff(B\hhh)} \ddd_{E,h/2},
$$
\end{itemize}
\end{defi}
The spaces  $\ddd_{B',h/2}$  and  $\ddd_{\hhh,h/2}$  are dense in  
$L^2(B,\mu _{B,h/2})$.
The functions  in $\ddd_{B',h/2}$  and  $\ddd_{\hhh,h/2}$ are defined on
$B$ but depend only on a finite number of variables. The slightly unnatural
parameter $h/2$
comes from \cite{AJN-JFA}, where the Segal-Bargman transformation links
$L^2(B,\mu_{B,h/2})$ and $L^2(B^2,\mu_{B,h})$. Such changes of variance are
unavoidable, see for example Prop. \ref{4.8AJN-JFA} and Def. \ref{13-AJNJFA}
below. 
\\

The Wigner functions, which are Gaussian Wigner
functions, are defined differently according to whether the
test functions are defined on  $E\in \fff(\hhh)$ or on $B$. The relationships
between these functions are specified in Prop. \ref{Wigner}.
The definition of $W_{h,E}$  (for test functions on $E$)
is inspired by the classical
definition of a Wigner function, taking into account the fact that the measure
is Gaussian.
\begin{defi}\label{Wigner-infini}
Let $E$ be in $\fff(\hhh)$, 
let  $\hat{f},\hat{g}$ be in  $\sss_{E,h/2}$.
The Wigner function of 
  $(\hat{f},\hat{g})$ is defined on $E^2$  by:
\begin{equation}\label{Wigner-gaussienne}
  \forall (z,\zeta)\in E^2,\quad 
W_{h,E}(\hat{f},\hat{g})(z,\zeta)=
e^{|\zeta|^2/h}\int_E e^{-2i\zeta \cdot t /h}
\hat{f}(z+t) \overline{\hat{g}(z-t)} e^{-|t|^2/h} (\pi h)^{-d/2} \ dt.
\end{equation}
Suppose that $ \tilde{f} $ and $\tilde{g}$ are defined on $B$ and satisfy 
$ \tilde{f}=\hat{f}\circ  \tilde{\pi}_E,  \tilde{g}= \hat{g}\circ  \tilde{\pi}_E$.\\
Then the   Wigner function of  $ \tilde{f}$ and
$ \tilde{g}$ is defined on $B^2$ by:
\begin{equation*}
\begin{array}{lll}
\displaystyle
\forall (z,\zeta)\in B^2,& \displaystyle
W_{h,B} ( \tilde{f}, \tilde{g}) (z,\zeta) &\displaystyle =
 W_{h,E}(\hat{f},\hat{g})(\tilde{\pi}_E(z),\tilde{\pi}_E(\zeta))\\ \\
 &  &\displaystyle = e^{|\tilde{\pi}_E(\zeta)|^2/2} 
\int_E e^{-2 \frac{i}{h} \tilde{\pi}_E(\zeta)\cdot t}\ \hat{f}(\tilde{\pi}_E(z)+t) 
\overline{\hat{g}(\tilde{\pi}_E(z)-t)} \ d\mu_{E,h/2}(t).
  \end{array}
\end{equation*}
\end{defi}
We need here the extended projections 
$\tilde{\pi}_E$. What we denote here by 
$
W_{h,E}(\hat{f},\hat{g}) $
would have been called
${H}_h^{Gauss}(\hat{f},\hat{g})$ in \cite{AJN-JFA},
formula (9), but we wish to
indicate on which space
the test functions are defined. \\

In the second case, the Wigner function  of  $(\tilde{f},\tilde{g})$ does
not depend on
the space $E$ on which  $\tilde{f}$ and $\tilde{g}$
are based.   We may state some of its properties, which hold for
$\tilde{f},\tilde{g}$ in  $\ddd_{\hhh,h/2}$ or $\ddd_{B',h/2}$. The result below
is taken from 
\cite{AJN-JFA} (Prop. 4.8). 
\begin{prop}\label{4.8AJN-JFA}
For all  $\tilde{f},\tilde{g}$ in $\ddd_{\hhh, h/2}$, the
Wigner function $W_{h,B}(\tilde{f},\tilde{g})$
belongs to $L^1(B^2, \mu_{B^2, h/2})$.
The operator associating, with all functions 
$\tilde{f},\tilde{g}$ of $\ddd_{\hhh,h/2}$, their Wigner function
$W_{h,B}(\tilde{f},\tilde{g})$,
extends uniquely as a continuous bilinear map from
$L^2(B,\mu_{B, h/2})\times L^2(B ,\mu_{B, h/2})$ in
 $L^2(B^2,\mu_{B^2, h/4})$, with norm $\leq 1$. 
\end{prop}

Now, we can give the definition of the quadratic form.
\begin{defi}\label{Weyl-quadratic-form}      
Let $\widetilde{F}$ be a bounded Borel function on $B^2$. One defines
$Q_h^{Weyl}(\widetilde{F}) $  by its action on $\ddd_{B',h/2}^2$: 
\begin{equation}\label{13-AJNJFA}
  \forall (\tilde{f},\tilde{g}) \in \ddd_{B',h/2}^2,\quad
Q_h^{Weyl}(\widetilde{F})  (\tilde{f},\tilde{g})
= \int _{B^2} \widetilde{F} (Z) W_{h,B}(\tilde{f},\tilde{g}) ( Z)
d\mu_{B^2,h/2} (Z) .
\end{equation}
If  $\widetilde{F}$  is not bounded, but if there exists an integer
 $m\geq 0$ such that
\begin{equation}\label{condition-N}
N_m(\widetilde{F}):=\sup_{Y\in H^2} 
\frac{||\tau_Y \widetilde{F}||_{L^1(B^2,\mu_{B^2,h/2})}}{(1+|Y|_{H^2})^m}
<+\infty, 
\end{equation}
then
$Q_h^{Weyl}(\widetilde{F}) $ can be defined as above. 
\end{defi}

Condition (\ref{condition-N}) means that the translates of  $\widetilde{F}$
by vectors $Y$  belonging to $\hhh^2$ still belong to $L^1$ and that their
norms depend polynomially at most on the translation vector $Y$.
If this vector were not in $\hhh^2$, the initial measure and the measure
obtained by translation would be mutually  orthogonal.\\

{\bf  Remarks:} 
This construction is a first step towards the construction of an operator.
But it is useful in itself, since it allows one to consider unbounded symbols
like polynomials in the functions $\ell_a$. For example, the symbol
 $\widetilde \varphi_{a,b} :$
$(x,\xi )\mapsto  \ell _a (x) + \ell _b (\xi) $
gives a multiplication by a monomial and a differentiation operator, as in the
finite dimensional case. \\
The normalization of the Wigner function $W_{h,E}$  has be chosen in order to
recover the classical Weyl calculus in the case when the symbols only depend
on a finite number of variables. Since we need a Gaussian measure, a
Gaussian factor appears in \eqref{Wigner-gaussienne}.\\

\noindent{\it Second construction }\\

A {\it multiindex} $\alpha$ is an element of  $\N^d$ when the dimension $d$ is
finite. Otherwise, it is a mapping from  $\N$  or $\N^*$ into $\N$ with
a finite number of nonzero coordinates. In both cases, one calls {\it depth}
of the multiindex the maximum of its coordinates, 
$\max_{j\in \{ 1,\dots, d\} }\alpha_j$ or $\max_{j\in \N^*}\alpha_j$.
\\

The first symbol class has features recalling the conditions of the 
C\`alderon-Vaillancourt Theorem and will be named after this result.
\begin{defi}\label{CV-Class}
 Let  $( \hhh,B)$  be an abstract Wiener space.
 Let $\bbb= (e_j )_{j\in \N^*}$ be a Hilbert basis of $\hhh$, each vector
 belonging to
$B'$. For $j\geq 1$ et $u_j = (e_j,0)$ and
$v_j = (0, e_j)$. Let $m$ be  a nonnegative integer
 and  $\varepsilon = (\varepsilon_j )_{j \in \N^*}$ be a family of
nonnegative real numbers.  One denotes by $ S_m(\bbb, \varepsilon)$ the
set of bounded continuous functions $ F:\hhh^2\rightarrow {\bf
  C}$ satisfying the following condition.
There exists $M\in \R^+$  such that, for any multiindices
 $\alpha,\beta$ of depth $m$, the following
derivative
\begin{equation}
\partial_u^{\alpha}\partial_v^{\beta}  F =  \left [\prod _{j\in \N^* }
\partial _{u_j} ^{\alpha_j} \partial _{v_j} ^{\beta_j}\right ]  F  \label{der-CV}
\end{equation}
is well defined, continuous on
 $\hhh^2$ and satisfies, for every $(x,\xi)$ in  $\hhh^2$
\begin{equation}
\left | \left [\prod _{j\in \N^* } \partial_{u_j} ^{\alpha_j}
 \partial _{v_j} ^{\beta_j}\right ]  F(x,\xi)
\right |    \leq M \prod _{j\in \N^* } \varepsilon_j ^{\alpha_j +
\beta_j}\ . \label{ineg-CV}
\end{equation}
Denote by $||F||= ||F||_{{ S_m(\bbb, \varepsilon)}}$ the smallest constant $M$
for which \eqref{ineg-CV} holds. With this norm
$|| \ ||_{ S_m(\bbb, \varepsilon)}$,
 $ S_m(\bbb, \varepsilon)$ is a Banach space. 
\end{defi}

One then has the existence and boundedness result 
  (\cite{AJN-JFA}, Th. 1.4).
\begin{theo}\label{1.4}
 Let  $( \hhh,B)$  be an abstract Wiener
space and let $h$ be a positive number. Let 
 $(e_j )_{j\in \N^*}$ be a Hilbert  basis of $\hhh$, each vector belonging
 to
$B'$. Let $F$  be a function on  $\hhh^2$
satisfying the following two hypotheses:
\begin{itemize}
\item
it belongs to the class $ S_2(\bbb, \varepsilon)$,
 where
  $\varepsilon =
(\varepsilon_j )_{(j \in \N^*)}$ is a square summable family of
nonnegative real numbers;
\item
it has a stochastic extension
$\widetilde F$ with   respect to both
 measures  $\mu_{B^2,h}$ and  $\mu_{B^2 ,h/2}$ (see Def.  \ref{ext-stoch}).
\end{itemize}
 Then there exists an operator,
denoted by  $Op_h^{Weyl}(F)$, bounded in  $L^2(B, \mu_{B, h/2})$,
such that, for all $\tilde{f}$ and $\tilde{g}$ in ${\cal D}_{B',h/2}$
\begin{equation}
  <Op_h^{Weyl}(F) \tilde{f},\tilde{g}> = Q_h^{Weyl}(\widetilde F)
  ( \tilde{f},\tilde{g}),\label{(1.17)} 
\end{equation}
where the right hand side is defined by Definition
 \ref{Weyl-quadratic-form}. Moreover,  if $h$ is in  $(0, 1]$:
\begin{equation} \Vert Op_h^{Weyl}(F)\Vert _{{\cal L}(
L^2(B,\mu_{B,h/2}) )}\leq ||F||_{ S_2(\bbb, \varepsilon)}  \prod_{j\in \N^*} (1 + 81 \pi h
S_{\varepsilon} \varepsilon _j^2),\label{(1.18)} 
\end{equation}
where
\begin{equation}
 S_{\varepsilon} = \sup _{j\in \N^*} \max (1 ,\varepsilon_j^2). \label{(1.19)}
\end{equation}
\end{theo}

We need the basis to belong to $B'$ because
of a decomposition result stated in Section \ref{sec-Gaarding}, hence the
use of the space ${\cal D}_{B',h/2}$.

In the more restrictive case when the sequence  $(\varepsilon_j) _{j\in \N^*}$
is summable, a function $F$ belonging to the C\`alderon-Vaillancourt class
$S_1(\bbb,\varepsilon)$ admits a stochastic extension in 
$L^q(B^2,\mu _{B^2,h})$
for all  $h>0$ and all  $q\in [1,+\infty[$ (\cite{J-StExt}, Prop. 3.1).
Moreover, there exists a function  $\tilde{F}$ which is the stochastic
extension of $F$ for all $h>0$ and all  $q \in [1,\infty[$.\\

Let us now define the second symbol class.
\begin{defi}\label{Classe-DD-def}
Let $A$ be a linear, selfadjoint, nonnegative, trace class application on a 
Hilbert space $\hhh^2$. For all $(x,\xi)\in \hhh^2$, one sets $Q_A(x,\xi)=
\langle A(x,\xi) , (x,\xi)\rangle $, where $\langle \ , \ \rangle$ denotes
the scalar product in $\hhh^2$.
Let  $S(Q_A)=S(Q_A,\hhh^2) $ be the class of all functions $F\in C^{\infty }
(\hhh^2)$ such that 
there exists $C(F) >0$ satisfying:
\begin{equation}\label{Classe-DD-def-ineq}
\begin{array}{lll}\displaystyle
 \forall (x,\xi)\in \hhh^2,\ 
 |F (x,\xi)  | \leq C (F),\\
 \displaystyle
\forall m\in\N^*, \forall (x,\xi)\in \hhh^2,  \forall (U_1,\dots,U_m)\in
(\hhh^{2})^m, 
|(d^m F ) (x,\xi) (U_1 , ... , U_m ) | \leq C (F) \prod _{j=1}^m
 Q_A( U_j) ^{\frac{1}{2}}.
\end{array}
\end{equation}
The smallest constant $C(F)$ such that  (\ref{Classe-DD-def-ineq})
holds is denoted by $\Vert F \Vert_{Q_A}$. 
\end{defi}

One checks that
$S(Q_A)$, endowed with the norm  $||\ ||_{Q_A}$, is a Banach space.\\

The class $S(Q_A)$ is more restrictive than
the preceding class.
Indeed, 
 $S(Q_A)\subset S_{\infty}(\bbb,\eps)$  for any orthonormal basis 
$\bbb= (e_j)$ of  $\hhh$, with
$\eps_j= \max( Q_A(e_j,0)^{1/2}, Q_A(0,e_j)^{1/2})$ and
$||F||_{S_2(\bbb,\eps)}\leq ||F||_{S_{\infty}(\bbb,\eps)}= \Vert F \Vert_{Q_A} $.
This sequence $\eps$ is only square summable. But
functions  in a class  $S(Q_A) $ admit a stochastic extension
in $L^p(B^2, \mu_{B^2, s})$ for all
$p\in [1,\infty[$ and $s>0$ for orthogonality
reasons  (see \cite{J-StExt}, Prop. 3.9). This justifies that  Theorem
\ref{1.4} holds for a symbol $F$ in  $S(Q_A)$.\\

\noindent{\bf Remark:}
 This class was introduced to relax an assumption in 
\cite{A-L-N-QED}. Moreover, its properties make the construction of the operator
easier, because one can use the Anti-Wick calculus - which exists in
the infinite dimensional frame - as a transition.  The results stated in
the next section are proved for the Calder\`on-Vaillancourt classes and
therefore valid for the classes of Definition \ref{Classe-DD-def}.

\subsection{Main results}\label{sec-results}

The first results of this section, Propositions
\ref{non-pos} and \ref{Flandrin1}, extend, to the Weyl calculus on the Wiener
space, 
results already known in the finite dimensional case. 
We show that an operator with a positive symbol is not necessarily positive,
exactly as in the case of Weyl's classical calculus.
The other result concerns Wigner functions and the 
Flandrin
conjecture.  Recalled in Part \ref{sec-Flandrin}, this conjecture
has been recently invalidated
in dimension $1$, for the configuration space $\R$,
in   \cite{DDL,Le}. Using these articles, we show that it does not hold either
when the configuration
space is the Wiener space.\\

The third result generalizes the paper  \cite{AJN-low}
concerning operators with a radial symbol,
in the finite dimensional case: if the symbol is radial, positive and
satisfies  further assumptions,
the operator is positive. This property holds for the
infinite dimensional calculus too.\\

The fourth result is probably the most satisfactory, since it is
 G\aa rding's inequality in the infinite dimensional frame. As in the classical
case, if the symbol $F$  is positive, we can't say that the operator
$ Op_h^{Weyl}(F)$ is positive. The loss of positivity is quantified in
Prop. \ref{Gaa}.

\begin{prop}\label{non-pos}
The Weyl calculus recalled in Part \ref{sec-Weyl-Wiener} is not positive.
There exists a positive symbol $F$ belonging to the 
Calder\`on -Vaillancourt classes of Definition \ref{CV-Class} and
a test function
$\tilde{g}\in L^2(B,\mu_{B,h/2})$ such that
$$
\langle  Op_h^{Weyl}(F)\tilde{g},\tilde{g}\rangle_{ L^2(B,\mu_{B,h/2})} <0.
$$
\end{prop}

Now the result about Flandrin's conjecture:
\begin{prop}\label{Flandrin1}
Take $e_1\in B'$ with norm $1$. For $a>0$ or  $a=+\infty$, define 
 $\tilde{F_a}$ on $B^2 $ by
$$
\tilde{F_a}(z,\zeta) = \1_{[0,a[}(\ell_{e_1}(\zeta))
\1_{[0,2\pi h a[}(\ell_{e_1}(z)) .
$$
For $a=\infty$, the indicator functions are $\1_{\R^+}$.\\ 
Then, for  $a>0$ sufficiently large or $a=\infty$, there exists a cylindrical function 
$\tilde{v_a}\in L^2(B,\mu_{B,h/2})$ such that
    
$$
Q_h^{Weyl}(\tilde{F_a})(\tilde{v_a},
\tilde{v_a})=
 \int_{B^2} \tilde{F}_a(z,\zeta) W_{B,h}(\tilde{v_a},
\tilde{v_a})(z,\zeta)
\ d\mu_{B^2,h/2}(z,\zeta) >
||\tilde{v_a}||^2_{L^2(B, \mu_{B,h/2})}    .
$$

\end{prop}
Observe that the function $F$ considered here does not belong to a symbol class,
since it is not defined on $\hhh^2$ and it is not even continuous.
Therefore, there is no operator associated with it. \\

For commodity reasons, the positivity result for symbols with radial properties
is stated below in an
expurgated form. A more general (and hence more technical) version,
Prop. \ref{prod-rad-gene},
is proved in
Part \ref{sec-posit-rad}.

\begin{prop}\label{posit-rad-gene}
Let $\Phi : \R^+ \rightarrow \R$  be a smooth, increasing, function.
Assume that
$\Phi$ is such that the radial function, defined on  $\R^{2d}$ by
$$
(x_1,\dots,x_d,\xi_1,\dots, \xi_d)\mapsto
\Phi(\sum_{j=1}^d(|x_j|^2 + |\xi_j|^2))$$
is smooth too. \\
Let $(e_1,\dots e_d)$ be an orthonormal family of $\hhh$, with all vectors
in $B'$. \\
Set
$$
\begin{array}{lll}
  \displaystyle
\forall (x,\xi)\in \hhh^2, 
F(x,\xi) = \Phi(\sum_{j=1}^d(|x\cdot e_j|^2 + |\xi \cdot e_j|^2)),\\ \\
\displaystyle
\forall(z,\zeta)\in B^2, \tilde{F}(z,\zeta)= 
\Phi(\sum_{j=1}^d(|\ell_{e_j}(x)|^2 + |\ell_{e_j}(\xi)|^2)).
\end{array}
$$
Theorem \ref{1.4} allows one to associate, with $F$, 
 an operator  $Op^{Weyl}_h(F)$, 
 bounded on  $L^2(B,\mu_{B,h/2})$. Moreover,
$Op^{Weyl}_h(F)$ has the following properties:
\begin{itemize}
\item
For all $\tilde{f} \in L^2(B,\mu_{B,h/2})$, one has
$$
\langle Op^{Weyl}_h(F) \tilde{f} ,\tilde{f}   \rangle_{L^2(B,\mu_{B,h/2})}    
\geq \left(  \frac{1}{h} \int_0^{\infty}
\Phi(t)e^{-t/h} \ dt\right)  \ ||\tilde{f} ||^2_{L^2(B,\mu_{B,h/2})}.
$$
\item
If, moreover, $\tilde{f} $ belongs to $\in \ddd_{B',h/2}$, the expression above
is given by the quadratic form:
$$
\langle Op^{Weyl}_h(F) \tilde{f} ,\tilde{f}   \rangle_{L^2(B,\mu_{B,h/2})}  =
\int_{B^2} \tilde{F}(z,\zeta) W_{h,B}(\tilde{f} ,\tilde{f} )(z,\zeta)
\  d\mu_{B^2, h/2}(z,\zeta) .
$$
\end{itemize}
Note that, under these assumptions, the positivity of the symbol $F$
implies the positivity of the operator. 
\end{prop}

We now state the  G\aa rding inequality, which holds for both symbol classes.
Its proof is not based  on a composition result, as will be seen in Section
\ref{sec-Gaarding}.\\

\begin{prop}\label{Gaa}
Let $\bbb=(e_j)_{j\geq 1}$ be an orthonormal basis of $\hhh$, with all vectors in $B'$.
Let $\eps=(\eps_j)_{j\geq 1}$ be a square summable
sequence of positive real numbers. Set
$S_{\eps}= \sup_{j\in \N^*} \max(1,\eps_j^2)$.\\
Let $F\in S_{2}(\bbb,\eps) $ be nonnegative on 
 $\hhh^2$. Suppose that
$F$ has a stochastic extension  $\tilde{F}$ for
 $L^2(B^2,\mu_{B^2, h/2})$ and
$L^2(B^2,\mu_{B^2, h})$ (which is the case if  $\eps$ is summable). 
\\
Then, for every function $\tilde{f} $ in
$L^2(B,\mu_{B, h/2})$, one has the following inequality
\begin{equation}\label{Gaa-ineq}
  \langle Op_h^{Weyl}({F})\tilde{f} ,\tilde{f}  \rangle_{L^2(B,\mu_{B, h/2})}
  \geq -||F||_{S_{2}(\bbb,\eps)} \sum_{j\geq 1}\lambda_j \
  \prod_{s\geq 1} (1+\lambda_s) \ ||\tilde{f} ||^2_{L^2(B,\mu_{B, h/2})},
\end{equation}
where $\lambda_j= 81\pi h S_{\eps} \eps_j^2$.
\end{prop}

Since the symbol classes defined by a quadratic form of Definition
\ref{Classe-DD-def}
are included in the  Calder\`on-Vaillancourt classes and since they have stochastic
extensions in $L^2$ for both measures, 
this G\aa rding
 inequality holds for them too: 

\begin{coro}\label{Gaa-DD}
Let $F$ belong to  $S(Q_A,\hhh^2) $ for a
linear, selfadjoint, nonnegative, trace class application $A$
defined on $\hhh^2$.\\
Then for all function $\tilde{f} $ in   $L^2(B,\mu_{B, h/2})$,
one has
\begin{equation}\label{Gaa-ineq-DD}
  \langle Op_h^{Weyl}({F})\tilde{f} ,\tilde{f}  \rangle_{L^2(B,\mu_{B, h/2})}
  \geq -  \Vert F \Vert_{Q_A}   \sum_{j\geq 1}\lambda_j \
  \prod_{s\geq 1} (1+\lambda_s) \ ||\tilde{f} ||^2_{L^2(B,\mu_{B, h/2})},
\end{equation}
where $\lambda_j= 81\pi h S_{\eps} \eps_j^2$ and
$S_{\eps}= \sup_{j\in \N^*} \max(1,   Q_A(e_j,0), Q_A(0,e_j))$.
\end{coro}


\section{Explicit stochastic extensions of cylindrical functions}\label{sec-Ext-stoch}

This section deals with stochastic extensions of cylindrical functions,
which are functions defined on $\hhh$ and depending on a finite number of
scalar products with elements of $\hhh$ (as recalled in Part \ref{sec-Wiener-space}). We prove that, in many cases,
one just needs to replace the scalar products by the corresponding
functions $\ell$  defined in  \eqref{application-ell}. This has been proved
for polynomial functions of scalar products in \cite{J-StExt}.\\
We first restate (for further reference) a lemma  which already appeared in
a previous article. It concerns the
functions $\ell$ themselves,
functions which replace
the monomials in the finite dimensional case.
We then give results
for cylindrical functions designed to play different parts in the calculus:
test
functions, symbols or Wigner functions.  \\

\begin{lemm}\label{cvprobella}
For every $a\in \hhh$, the  scalar product function $\hhh\rightarrow \R, x\mapsto a\cdot x$
admits, as a stochastic extension in  $L^p(B,\mu_{B,s})$, the  function
$\ell_a$. This holds for all $p\in [1,+\infty[$ and $s>0$. \\
This means that, if $(E_n)_{n\in \N}$ in an increasing sequence of $\fff(\hhh)$,
with union dense in  $\hhh$, the sequence of random variables 
$ (\ell_a - \ell_{P_{E_n}(a)})_{n\in \N}$ converges to $0$
in $L^p(B,\mu_{B,s})$ (for all  $p\in [1,+\infty[$ and $s>0$).
Hence it converges in 
$\mu_{B,s}$ probability too.
\end{lemm}

\noindent{\bf Proof}\\
For $E\in \fff(\hhh)$, one checks that
\begin{equation}\label{proj}
a \cdot \widetilde{\pi_E}= \ell_{P_E(a)}.
\end{equation}
Indeed, if $(e_1\dots, e_d)$ is an orthonormal basis of $E$ for the scalar
product of $\hhh$,
$$
\forall x\in B,\quad a \cdot \widetilde{\pi_E}(x)= a\cdot \sum_1^d \ell_{e_i}(x) e_i =
\sum_1^d \ell_{e_i}(x) a\cdot e_i =
\ell_{ \sum_1^d (a\cdot e_i) e_i }(x)
$$
by linearity of $a \mapsto \ell_a$. One then recognizes the orthogonal
projection on $E$.\\
Let $(E_n)_{n\in \N}$  be an increasing sequence of  $\fff(\hhh)$,
with union dense in
$\hhh$. 
For $1\leq p< +\infty$, the formula above and Definition \ref{ext-stoch}
of stochastic extensions lead us to consider
$||\ell_a - \ell_{P_{E_n}(a)}||_{L^p(B,\mu_{B,s})} $. We may write that
$$
||\ell_a - \ell_{P_{E_n}(a)}||_{L^p(B,\mu_{B,s})} = C_{p,s} |a-P_{E_n}(a)|,
\quad  {\rm with  }\
C_{p,s}=\sqrt{2s}\ \pi^{-1/2p}\ 
\Gamma(\frac{p+1}{2})^{1/p}.
$$
This is a consequence of the transfer theorem, which brings us back to
the following integral on $\R$:
\begin{equation}\label{ell-a-Lp}
  \forall b\in \hhh\setminus \{0\}, \quad 
  \int_B|\ell_b|^p d\mu_{B,s}=
  \int_{\R} |x|^p  (2\pi s |b|^2)^{-1/2} e^{-\frac{x^2}{2s|b|^2}} dx =(2s)^{p/2}\ \pi^{-1/2}\ \Gamma(\frac{p+1}{2})\ |b|^p.
\end{equation}
Hence
$||\ell_a - \ell_{P_{E_n}(a)}||_{L^p(B,\mu_{B,s})} $
converges to  $0$, which, in turn, implies the convergence in probability.
\hfill $\square$\\

We now turn to  more general regular cylindrical functions on $\hhh^2$.
Under regularity
and decay assumptions, they belong to a 
C\`alderon-Vaillancourt symbol class (of Definition \ref{CV-Class}).
Even if this gives the existence of the stochastic extension when the
sequence $\eps$ is summable, 
the extension is not necessarily explicit. In the case of cylindrical
functions, one may  be more precise. 

\begin{lemm}\label{dans-SBeps}
Let $\Ddot{F}$ be a bounded,   $C^{m}$ function  on $\R^{2d}$, with
bounded partial derivatives of all orders (smaller than $m$).
Let  $(e_n)_{n\in \N^*}$
be a Hilbert basis of $\hhh$, with elements in $B'$.\\
Let  the functions $F$ and  $\tilde{F}$  be defined, respectively,
on  $\hhh^2$ and  $B^2$ by:
\begin{equation*}
  \begin{array}{lll}
     \displaystyle
\forall (x,\xi)\in \hhh^2, \quad {F} (x,\xi)=
\Ddot{F}( e_1 \cdot  x ,\dots,
 e_d\cdot x,
 e_1\cdot \xi ,\dots,
 e_d\cdot \xi),\\
    \displaystyle
\forall (z,\zeta)\in B^2, \quad \tilde{F}(z,\zeta)=
\Ddot{F}(\ell_{e_1}(z),\dots, \ell_{e_d}(z),\ell_{e_1}(\zeta),\cdots,
\ell_{e_d}(\zeta)).\\
  \end{array}
\end{equation*}
Then 
$F$ belongs to the symbol class  $ S_m(\bbb,\eps)$ for the sequence  
$\eps= (n^{-2})_{n\geq 1}$
\\
The function $F$ admits  $\tilde{F}$ as a stochastic extension in
$L^p(B^2,\mu_{B^2,s})$ for all finite $p\geq 1$ and all $s>0$. 
\end{lemm}

\noindent{\bf Proof }\\
Let $\alpha$ and $\beta$ be two multiindices of depth $m$.
If  $\alpha$ or $\beta$  has a nontrivial component for
an index $j>d$, inequality \eqref{ineg-CV} holds because its
left side is equal to $0$. Otherwise, there is a finite number of inequalities
to satisfy and it suffices to set
$$
M= \sup_{\alpha,\beta} \left(  \prod_{j\leq d}j^{2(\alpha_j+\beta_j)} \sup_{(x,\xi)\in \R^{2d}} \left|
\prod_{j=1}^d \frac{\partial^{\alpha_j}}{\partial x_j^{\alpha_j}}
\frac{\partial^{\beta_j}}{\partial \xi_j^{\beta_j}}\Ddot{F}(x,\xi)\right|\right) .
$$
The supremum above is taken on all multiindices of depth $m$ with
all components equal to $0$ for $j\geq d$.\\
Since the sequence  $\eps$ is summable, $F$ has a stochastic extension 
(temporarily denoted by  $F^*$)   in $L^p(B^2,\mu_{B^2,s})$ for all 
$p\in [1,\infty[$ and all $s>$, according to Proposition 3.1 of \cite{J-StExt}.
We will prove that it coincides  with  $\tilde{F}$.\\    
Let $(E_n)_{n\in \N^*} $ be an increasing sequence of $\fff(\hhh)$, with union dense in
$\hhh$. 
Set $F_n(z,\zeta)=F(\tilde{\pi}_{E_n}(z), \tilde{\pi}_{E_n}(\zeta))$, for $z,\zeta\in B$. 
Since   $(F_n)_{n\in \N^*}$ converges to  $F^*$ in
$L^p(B^2,\mu_{B^2,s})$,  a subsequence 
$(F_{\ph(n)})$ converges $\mu_{B^2,s}$-almost surely to  $F^*$.
According to the definition of  $F$, one has
$$
F_{\ph(n)}(z,\zeta)=
\Ddot{F}( e_1\cdot \tilde{\pi}_{E_{\ph(n)}}(z) ,\dots,
e_d\cdot \tilde{\pi}_{E_{\ph(n)}}(z),
e_1 \cdot \tilde{\pi}_{E_{\ph(n)}}(\zeta) ,\dots,
e_d \cdot \tilde{\pi}_{E_{\ph(n)}}(\zeta) )
$$
Lemma  \ref{cvprobella} gives the convergence in  $\mu_{B,s}$  probability of 
$ e_j\cdot \tilde{\pi}_{E_{\ph(n)}}(\cdot ) = \ell_{P_{E_{\ph(n)}}(e_j)}$
to $\ell_{e_j}$. Extracting a further subsequence and using the continuity of
$\Ddot{F}$, we get that
$$
F_{\ph(\psi(n))}(z,\zeta)\ \longrightarrow \
\Ddot{F}(\ell_{e_1}(z),\dots, \ell_{e_d}(z),\ell_{e_1}(\zeta),\cdots,
\ell_{e_d}(\zeta)) =\tilde{F}(z,\zeta) 
$$
$\mu_{B^2,s}$-almost surely. Hence, the  functions $F^*$ and $\tilde{F}$
are almost surely equal  for every measure  $\mu_{B^2,s}$,  which proves 
Lemma \ref{dans-SBeps}. 
\hfill $\square$\\

We now prove that the  functions of  $\ddd_{\hhh, h/2}$, defined in Definition
\ref{space-D}  are stochastic extensions of functions cylindrical on $\hhh$.
But we may state a slightly more general result, namely:
\begin{theo}\label{ext-stoch-1}
Let $s>0$.
Let $p\geq 1$.
Let $\ph : \R^d \rightarrow \R$ be  continuous and such that the function
\begin{equation}\label{conditionp}
  x \mapsto e^{-|x|^2/2ps} \ph(x)  P(x)
\end{equation}
is bounded for every polynomial $P$. \\
Let $E\in \fff(\hhh)$ have dimension $d$ and an orthonormal basis
$(e_1,\dots, e_d)$.\\
Define $f$ on  $\hhh$ and $\tilde{f}$ on  $B$  setting:
$$\forall x \in \hhh,\quad f(x)= \ph( (x\cdot e_i)_{1\leq i\leq d})\quad
{\rm and} \quad
\forall y\in B, \quad  \tilde{f}(y) = \ph( (\ell_{e_i}(y))_{1\leq i\leq d}).
$$
Then $\tilde{f}$ is  the stochastic extension  of $f$ in
$L^p(B,\mu_{B,s})$.\\
\end{theo}

This result has the following important corollaries :
\begin{coro}\label{enfin} 
Let $E\in \fff(\hhh)$, let  $s>0$. Let 
$(e_1,\dots,e_d)$ be an orthonormal basis of  $E$. If
$f$, defined on $\hhh$, has the form 
$ f(x)= \ph( (x\cdot e_i)_{1\leq i\leq d})$ with  $\gamma_{\R^d,s}\ph$
rapidly decreasing (which means that $f\in \sss_{E,s}$),
then $f$ has a stochastic extension in 
 $L^2(B,\mu_{B,s})$, which is the   function $\tilde{f}$ defined on $B$
by
$ \tilde{f}(y) = \ph( (\ell_{e_i}(y))_{1\leq i\leq d})$.
\end{coro}

This is true because, if $f$ satisfies the above conditions, 
$x\mapsto e^{-|x|^2/4s}\ph(x) $ is rapidly decreasing.
The exponent $p$ of Theorem \ref{ext-stoch-1}  is equal to $2$, the 
transformation $\gamma$ is defined in \eqref{gamma}. \\
This proves that the functions belonging to $\ddd_{\hhh,s}$ are indeed
stochastic extensions of functions defined thanks to the same $\ph$,  scalar
products with $a\in \hhh$  replacing the functions $\ell_a$.\\

\begin{coro}
Let $\ph : \R^d \rightarrow \R$ be continuous, compactly supported,
let $E\in \fff(\hhh)$ have dimension $d$ and an orthonormal basis
$(e_1,\dots, e_d)$.\\
Define $f$ on  $\hhh$ and $\tilde{f}$ on  $B$ setting:
$$\forall x \in \hhh,\quad f(x)= \ph( (x\cdot e_i)_{1\leq i\leq d}),\qquad
\forall y\in B, \quad  \tilde{f}(y) = \ph( (\ell_{e_i}(y))_{1\leq i\leq d}).
$$
Then, for all 
 $p\geq 1$ and  $s>0$, 
$\tilde{f}$ is  the stochastic extension  of $f$ in 
$L^p(B,\mu_{B,s})$.
\end{coro}
This result holds because, if $\ph$ is continuous and compactly supported,
the function defined in 
\eqref{conditionp} for a polynomial $P$
is bounded for all $p\in [1,+\infty[$ and $s>0$.\\
  
\noindent{\bf Proof  of Theorem  \ref{ext-stoch-1}}\\
The  $d$-uple $(\ell_{e_1},\dots,\ell_{e_d})$ has distribution
$\nnn(0, s I_d)$. Hence, 
$$
\int_B|\tilde{f}(x)|^p \ d\mu_{B,s}(x) =
\int_{\R^d} |\ph(y)|^p (2\pi s)^{-d/2} e^{-|y|^2/2s} \ dy.
$$
Condition \eqref{conditionp} on $\ph$ ensures that
the integral is finite, which shows that  $\tilde{f}\in L^p(B,\mu_{B,s})$.
\\
Recall that, in  a measured space 
 $(X,\ttt,\mu)$, 
if  $p\in [1,\infty[$ and if $(f_n)_{n\in \N }$ is a sequence of $\lll^p(\mu)$,
    converging almost everywhere to  $f\in \lll^p(\mu)$, then
\begin{equation}\label{BrPa}
\lim_{n\rightarrow \infty}||f-f_n||_{p} = 0 \quad \Longleftrightarrow \quad
\lim_{n\rightarrow \infty}||f_n||_{p} =||f||_{p} .
\end{equation}
(See, for example,  \cite{BP}).
Let us then take an increasing sequence of 
$\fff(\hhh)$, $(E_n)_{n\in \N }$ with union dense in $\hhh$.
Formula \eqref{proj}  implies that, $\mu_{B,s}$-almost everywhere on $B$,
\begin{equation}\label{P-rond-pitilde}
P_E(\tilde{\pi}_{E_n}(x))= \sum_{i=1}^{d} \ell_{P_{E_n}(e_i) }(x) e_i,
\end{equation}
where $P_E$ is the orthogonal projection on $E$, defined by \eqref{projection}. 

According to Lemma \ref{cvprobella}, for $1\leq i\leq d$, 
$\ell_{P_{E_n}(e_i) } $ converges to  $\ell_{e_i}$ in
$\mu_{B,s}$ probability.

Take a subsequence of  $(E_n)$ indexed by $\psi(n)$. There exists a further subsequence, indexed by  $\psi(\zeta(n))$, such that, for all $i\leq d$, 
$\ell_{e_i} - \ell_{P_{E_{\psi(\zeta(n)) }}(e_i)}$ converges $\mu_{B,s}$-almost everywhere
to $0$.
By the definition of $f$ and by formula \eqref{proj}, for every $y\in B$,
$$
f(\tilde{\pi}_{E_{\psi(\zeta(n))}}(y))
=\ph((\tilde{\pi}_{E_{\psi(\zeta(n))}}(y) \cdot e_i)_{i\leq d})=
\ph( ( \ell_{P_{E_{\psi(\zeta(n)) }}(e_i)}(y))_{i\leq d}).
$$
Since $\ph$ is  continuous, we get that
$$
f\circ {\tilde{\pi}}_{E_{\psi(\zeta(n))}} \quad \longrightarrow \quad
\ph((\ell_{e_i})_{i\leq d}) = \tilde{f}\quad \mu_{B,s} \ p.s.
$$
We now must check that
$|| \ph( ( \ell_{P_{E_{\psi(\zeta(n)) }}(e_i)})_{i\leq d}) ||_{L^p(B,\mu_{B,s})}$ converges
to $||\tilde{f}||_{L^p(B,\mu_{B,s})}$.

The $d$-uple $(\ell_{P_{E_{\psi(\zeta(n)) }}(e_i)})_{i\leq d})$ is normally distributed,
with  $0$ means and covariance matrix   $K_n$ equal to 
\begin{equation*}\label{Kn}
K_n= (( \langle \ell_{P_{E_{\psi(\zeta(n))}(e_i)}},  \ell_{P_{E_{\psi(\zeta(n))}(e_j)}}\rangle_{L^2(B,\mu_{B,s})})_{i,j\leq d} =
 s ( P_{E_{\psi(\zeta(n))}(e_i)}\cdot  P_{E_{\psi(\zeta(n))}(e_j)} )_{i,j\leq d},
\end{equation*}
where $\cdot$ is is the scalar product of    $\hhh$.
When $n$ goes to infinity, the coordinates of $K_n$  converge to 
$s \ e_i\cdot e_j$ and  $K_n$ itself converges to $s I_d$.
For sufficiently large $n$, $K_n$ is then invertible, hence
$( \ell_{P_{E_{\psi(\zeta(n)) }}(e_i)})_{i\leq d}$  admits a
density and 
$$
||f\circ  {\tilde{\pi}}_{E_{\psi(\zeta(n))}} ||^p_{L^p(B,\mu_{B,s})}=
\int_{B} | \ph( ( \ell_{P_{E_{\psi(\zeta(n)) }}(e_i)})_{i\leq d}) |^p \ d \mu_{B,s}=
\int_{\R^d} |\ph(y)|^p \frac{1}{(2\pi)^{d/2} \sqrt{{\rm det}(K_n)}}
e^{-\frac{1}{2} <y, K_n^{-1} y>}  \ dy.
$$
We know that, for all  $y\in \R^d$, 
$$
s <y, K_n^{-1} y>\ \geq \  <y,y>= |y|^2
$$
with the scalar product and the Euclidean norm on  $\R^d$.
This fact will be proved in Lemma  \ref{varcovar} below.
Moreover, since ${\rm det}(K_n)$ converges to $s^d$, its inverse is
bounded independently of $n$ for $n $ large enough. All this  allows
using the dominated convergence Theorem, since  $K_n^{-1}$ converges to
$s^{-1} I_d$. One deduces that
$$
||f\circ  {\tilde{\pi}}_{E_{\psi(\zeta(n))}} ||^p_{L^p(B,\mu_{B,s})}\quad
\longrightarrow \quad
\int_{\R^d} |\ph(y)|^p (2\pi s)^{-d/2} 
e^{-\frac{1}{2s} <y, y>}  \ dy =
||\tilde{f}  ||^p_{L^p(B,\mu_{B,s})}.
$$
Then, for this subsequence, \eqref{BrPa} yields that
$$
||f\circ {\tilde{\pi}}_{E_{\psi(\zeta(n))}} -\tilde{f} ||_{L^p(B,\mu_{B,s})}\quad
\longrightarrow \quad 0.
$$
A proof by contradiction  then ensures that
$
||f\circ  {\tilde{\pi}}_{E_n} -\tilde{f} ||_{L^p(B,\mu_{B,s})}$
itself converges to $0$, which achieves the proof of
Theorem  \ref{ext-stoch-1}.
\hfill $\square$\\

We now give the result about the inverse of $K_n$.
\begin{lemm}\label{varcovar} 
Let $K_n=  s ( P_{E_{n}(e_i)}\cdot  P_{E_{n}(e_j)} )_{i,j\leq d}
$ be the covariance matrix appearing in the preceding proof. 
For all $y\in \R^d$,
$$
s <y, K_n^{-1} y>\  \geq \ <y,y>= |y|^2.
$$
\end{lemm}
\noindent{\bf Proof  :}\\
For $n$ large enough, $K_n$ is invertible. According to its definition, it is definite positive and its inverse has the same property. Using  their
square roots gives
$
<y, K_n^{-1} y>=  < K_n^{-1/2}y, K_n^{-1/2} y>.
$
Setting 
 $x=  K_n^{-1/2} y$ we get
$$
\begin{array}{lll}
\forall y\in \R^d, \ s  \langle y, K_n^{-1} y\rangle \geq  \langle y,y\rangle  &
\Longleftrightarrow &\displaystyle 
\forall x \in \R^d, \ \langle x,x\rangle \ \geq \  \frac{1}{s}
\langle  K_n^{1/2}x, K_n^{1/2}x\rangle\\ \\
&
\Longleftrightarrow &\displaystyle 
\forall x \in \R^d, \  \langle x,x\rangle \geq \frac{1}{s} \langle x, K_n x\rangle. \\
\end{array}
$$
Since
$$
\langle x, K_n x\rangle = s\sum_{i,j} x_ix_j  P_{E_n}(e_i)\cdot P_{E_n}(e_j)
= s\sum_i x_i P_{E_n}(e_i)\cdot \sum_j x_j P_{E_n}(e_j),
$$
one has
$$
\langle x, K_n x\rangle = s|P_{E_n}(\sum_i x_i e_i)|_{\hhh}^2
\leq s |\sum_i x_i e_i|_{\hhh}^2 = s\langle x,x\rangle.
$$
This proves the inequality. 
\hfill $\square$\\

We now apply the extension results to the Wigner functions.
Let $E\in \fff(\hhh)$.
According to Def. \ref{Wigner-infini}, for $\hat{f},\hat{g}\in \sss_{E,h/2}$,
$W_{h,E}(\hat{f},\hat{g})$ is given by
 $$
 \forall (z,\zeta)\in E^2,\quad
W_{h,E}(\hat{f},\hat{g})(z,\zeta)=
e^{|\zeta|^2/h}\int_E e^{-2i\zeta \cdot t /h}
\hat{f}(z+t) \overline{\hat{g}(z-t)} e^{-|t|^2/h} (\pi h)^{-d/2} \ dt.
$$
Define functions $f,g$ on $\hhh$ by $f=\hat{f}\circ P_E,g=\hat{g}\circ P_E$,
let $\tilde{f},\tilde{g}$ be defined on $B$ by
$\tilde{f}= \hat{f}\circ \tilde{\pi}_E, \tilde{g}= \hat{g}\circ \tilde{\pi}_E$.
Corollary \ref{enfin} says  that $\tilde{f}$
is the stochastic extension of $f$  in $L^2(B,\mu_{B,h/2})$.\\

According to Def. \ref{Wigner-infini} again,
the function $W_{h,B}(\tilde{f},\tilde{g})$
is given by 
$$\forall (z,\zeta)\in B^2,\quad
W_{h,B}(\tilde{f},\tilde{g})(z,\zeta)=
W_{h,E}(\hat{f},\hat{g})(\tilde{\pi}_E(z), \tilde{\pi}_E(\zeta)) .$$
We now set 
$$
\forall (z,\zeta)\in \hhh^2,\quad
W_{h,\hhh}(f,g)(z,\zeta) = W_{h,E}(\hat{f},\hat{g})(P_E(z),P_E(\zeta))
.
$$

One sees that, to get $W_{h,B}(\tilde{f},\tilde{g})$ when one knows
$W_{h,\hhh}(f,g)$,
it suffices to replace the scalar products (in the  projections $P_E$)
by the corresponding functions $\ell$.\\

We then have the following result. 
\begin{coro}\label{Wigner}
 With the preceding notations, 
 $W_{h,B}(\tilde{f},\tilde{g})$ is  the stochastic extension
 of $W_{h,\hhh}(f,g)$  in
$L^2(B^2,\mu_{B,h/4})$.
\end{coro}

\noindent{\bf Proof of Corollary \ref{Wigner}}\\
We express the Wigner function thanks to
$\gamma_{E,h/2}\hat{f}$ and $\gamma_{E,h/2}\hat{g}$,
which are rapidly decreasing. This yields 
$$
W_{h,E}(\hat{f},\hat{g})(z,\zeta)=
e^{(|z|^2+|\zeta|^2)/h}\int_E e^{-2i\zeta \cdot t /h}
\gamma_{E,h/2}\hat{f}(z+t) \overline{\gamma_{E,h/2}\hat{g}(z-t)} \ dt.
$$
The integral factor above is rapidly decreasing in  $(z,\zeta)$ (see,
for example, \cite{DDL} for a recent reference). Applying Corollary  \ref{enfin}
 with  $s=h/4$ for the ``target''  space $B^2$ proves that 
 $W_{h,\hhh}(f,g)$ has a stochastic extension   in $L^2(B^2, \mu_{B^2, h/4})$,
 which is $W_{h,B}(\tilde{f},\tilde{g})$. \\


 \section{Around  positivity }\label{sec-pos}

This section extends, to the infinite dimension, positivity or nonpositivity
results already known in the finite dimensional case.
The main tools are the finite dimensional Hermite functions,
their Wigner functions, recalled
in the first part below, their stochastic extensions, studied in the
preceding Section \ref{sec-Ext-stoch}.\\
The point is that results in finite dimension  can be
transposed to cylindrical functions (resp. cylindrical symbols),
officially defined on
an infinite dimensional space, but in fact depending on $d$
(resp. $2d$) variables.

\subsection{The Wigner functions of the  Hermite functions - finite dimensional case}
\label{basedim1}

We first fix the normalization for the Hermite functions in dimension $1$.
We give a Hilbert basis of $L^2(\R, \mu_{\R,h/2})$, to make the
transition with the choices made in the infinite dimensional case easier.
The definition of the Hermite functions and some of their properties are
followed by
the computation of their Wigner functions. The formulas are given explicitly,
for sake of clarity. Since they are classical (up to normalization) and may be
found, with their proofs,  in \cite{F,Jan}, the arguments are
only sketched here. At the end of this part, a technical lemma allowing
integrations by parts is  proved. Il will be applied in Part \ref{sec-posit-rad}
about radial symbols. 
\\

One denotes by $\mu_{\R,h/2}$ the Gaussian measure with density
$(\pi h)^{-1/2} e^{-x^2/h}$ with respect to the Lebesgue measure. One sets
$\tau(x)=  e^{-\frac{x^2}{h}}$. For every real number $x$, one sets 
\begin{equation}\label{psi}
  \psi_{-1}(x)=0, \quad
\forall j\geq 0,\ 
\psi_j(x)= \frac{ (-1)^j}{\sqrt{ j!}}
\left(\frac{h}{2}\right)^{j/2}\  e^{\frac{x^2}{h}} \ \tau^{(j)}(x).
\end{equation}
The classical relations then
have this form:
\begin{equation}\label{eq-rels}
\begin{array}{llllll}
   \forall x \in  \R,\quad &
\displaystyle \forall j\geq 0,&  \displaystyle 
\psi_{j+1}(x)=  \frac{-1}{\sqrt{j+1}}\sqrt{\frac{h}{2}}\left(\left(
\frac{d}{dx} - \frac{2}{h} x\right)\psi_j\right)  (x) & {\rm creation}
\\ \\
&\displaystyle \forall j\geq 1,& 
\displaystyle
\psi_j(x)= \sqrt{\frac{2}{h}}\frac{x}{\sqrt{j}} \ \psi_{j-1}(x)
-\sqrt{\frac{j-1}{j}}\  \psi_{j-2}(x)  & {\rm recurrence}\\
\\
& \displaystyle \forall j\geq 1,& 
\displaystyle
\frac{d \psi_j}{dx} = \sqrt{j}  \sqrt{\frac{2}{h}}\ \psi_{j-1}
& {\rm annihilation}
\end{array}
\end{equation}
With $\psi_{-1}(x)=0$, the first Hermite functions are, explicitly,
\begin{equation}\label{Herm-prem}
\begin{array}{lll}
\displaystyle
\psi_0(x)=1  &   \displaystyle
\psi_1(x)= \sqrt{\frac{2}{h}} \ x \\ \\
\displaystyle
\psi_2(x)=\frac{1}{\sqrt{2}} \left( \frac{2}{h} x^2-1\right)  &
\displaystyle
\psi_3(x)=\frac{1}{\sqrt{6}}  \left( \left(\frac{2}{h}\right)^{3/2} x^3
-3\left(\frac{2}{h}\right)^{1/2}x \right)  \\ \\
\end{array}
\end{equation}

The recurrence formula proves that
$\psi_j$ is a polynomial with degree $j$ exactly and leading coefficient 
$\displaystyle \left( \frac{2}{h} \right)^{j/2} \frac{1}{\sqrt{j!}}.$
The family  $(\psi_j)_{j\geq 0}$  is an orthonormal family of
$L^2(\R, \mu_{\R, h/2})$.\\


\noindent{\it Computation of the Wigner functions of the Hermite functions in dimension 1}\\

Direct computations and  the definition formula
 \eqref{Wigner-gaussienne} for $E=\R$ lead to the following results
\begin{equation}\label{Wig-prem}
 W_{h,\R}(\psi_0,\psi_0)(x,\xi) = 1,
 W_{h,\R}(\psi_0,\psi_1)(x,\xi) =\sqrt{\frac{2}{h}} (x+i\xi),
 W_{h,\R}(\psi_1,\psi_1)(x,\xi) =-1 + \frac{2}{h} (x^2+\xi^2).
\end{equation}
For more general orders, one needs to introduce the {\it Bargman
  kernel}  $K_v$ associated  with a complex number $v$ and defined
below:
for every  $v\in \C$ and every $x\in \R$, set
\begin{equation}\label{noy-B}
K_v(x) =\sum_{j=0}^{\infty}\psi_j(x) \frac{v^j}{\sqrt{j!}} = e^{xv\sqrt{2/h} - v^2/2}.
\end{equation}
The second equality comes from Taylor's Formula. 
For a fixed  $v\in \C$, the convergence takes place in
$L^2(\R,\mu_{\R, h/2})$ because the sequence of the coefficients
of the  $\psi_j$ is square summable. For a fixed real number $x$,
the convergence is uniform in every compact subset of $\C$, because the
convergence radius (in the variable $v$) is infinite.
\\
One checks that the following Wigner function of two Bargman kernels
is equal to
\begin{equation}\label{Wigner-noy-B}
W_{h,\R} (K_u, K_{\bar v})(x,\xi) =
\exp\big( -uv + \sqrt{\frac{2}{h}}x(u+v) +i \sqrt{\frac{2}{h}}\xi(v-u)\big).
\end{equation}

Exchanging discrete sums and integration and identifying the coefficient  of
the term $u^jv^k$, one proves that, for all real numbers   $x,\xi$
and all integers $j,k$:
\begin{equation}
W_{h,\R}(\psi_j,\psi_k)(x,\xi) =
\sum_{q= \max(0,k-j)}^k  \frac{(-1)^{q-k}}{(k-q)! }
\left( \frac{2}{h}\right)^{(j-k+2q)/2}
\frac{1}{(j-k+q)! q! } (x-i\xi)^{j-k+q} (x+i\xi)^q \sqrt{j!k!}
\end{equation}
Hence
\begin{equation}\label{Wig-dim1}
\begin{array}{lll} 
\displaystyle W_{h,\R}(\psi_j,\psi_k)(x,\xi) =
\sqrt{\frac{j!}{k!}}(x+i\xi)^{k-j}(-1)^j 
\left( \frac{2}{h}\right)^{(k-j)/2}
L_j^{(k-j)}\left(  \frac{2}{h}(x^2+\xi^2)\right) & {\rm if} & j\leq k\\ \\
\displaystyle W_{h,\R}(\psi_j,\psi_k)(x,\xi) =
\sqrt{\frac{k!}{j!}}(x-i\xi)^{j-k}(-1)^k 
\left( \frac{2}{h}\right)^{(j-k)/2}
L_k^{(j-k)}\left(  \frac{2}{h}(x^2+\xi^2)\right) & {\rm if} & j\geq  k.\\ \\
\end{array}
\end{equation}

Up to normalization, these results and their proofs can be found in 
 \cite{F}, Theorem 1.105. 
 Laguerre polynomials are defined in 
\cite{F, MOS} by: 
\begin{equation}\label{Laguerre}
L_k^{(\alpha)}(x) =
\sum_{m=0}^k \frac{(k+\alpha)! }{(k-m)! (\alpha +m)! }  \frac{(-x)^m}{m!}
\end{equation}

The leading coefficient of $L_k^{(\alpha)}$  in \eqref{Laguerre}
is  $\frac{(-1)^k}{k!}$. The
Laguerre polynomials $(L_k^{(\alpha)})_k$  are orthogonal with respect to
the measure
$ x^{\alpha} e^{-x}\ dx $ over $\R^+$, but not orthonormal.
They are normalized by their leading coefficient. \\

The last result in dimension $1$ is an integration by parts lemma, which
applies to the radial functions of a section below, or more generally to
functions cancelled by the differential operator
$[x\partial_{\xi} -\xi\partial_x]$.

\begin{lemm}\label{IPP}
Let $P : \R\rightarrow \R$ be a $C^1$   function, such that  $P $ and
$P'$ are at most polynomially increasing.\\
Let $s\in \N$, let $n\in \N^*$. Let $\eps=1 $ or $-1$. \\
Let $F :  \R^2\rightarrow \R$  be a $C^{n}$ function,
at most polynomially increasing, along  with its partial derivatives.\\
One has:
\begin{equation}
(-s i \eps)^n \int_{\R^2} F(x,\xi) (x+i\eps \xi)^s P(x^2 + \xi^2) \
e^{-\frac{x^2+\xi^2}{h}} \ dx d\xi =
\int_{\R^2}([x\partial_{\xi} -\xi\partial_x]^n F)(x,\xi) (x+i\eps \xi)^s P(x^2 + \xi^2) \
e^{-\frac{x^2+\xi^2}{h}} \ dx d\xi .
\end{equation}

\end{lemm}
\noindent{\bf Proof }\\ We start from the right hand-side.\\
In the case when $s\neq 0$, the properties of $F,P$
and the presence of the exponential function allow one to integrate by parts. 
The operator   $[x\partial_{\xi} -\xi\partial_x]$ cancels radial functions.
Hence, 
the only  term that does not vanish comes from
$-[x\partial_{\xi} -\xi\partial_x] (x+i\eps \xi)^s =(-i\eps s)(x+i\eps \xi)^s $.
In particular, no derivative of $P$ remains.
\\
If $s=0$, the equality has the form $0=0$. Indeed, $s=0$ is a factor of the
left term and an integration by parts cancels the right term, since
$[x\partial_{\xi} -\xi\partial_x]$ applies only to radial functions. 
\hfill $\square$\\

We can use this lemma  and the expression of $ W(\psi_j,\psi_k)$
to prove the following equalities, which  are the classical equalities
in the finite dimensional case:
\begin{equation}\label{Wigner-delta}
  \forall (j,k)\in \N^2,\quad
\int_{\R^2} W(\psi_j,\psi_k) (x,\xi) \ d\mu_{\R^2,h/2}(x,\xi) = \delta_{j,k}.
\end{equation}
Indeed,
if $j<k$, one has
$$
\int_{\R^2} W(\psi_j,\psi_k) (x,\xi) \ d\mu_{\R^2,h/2}(x,\xi)
=
C \int_{\R^2} (x+i\xi)^{k-j} L_j^{(k-j)}( \frac{2}{h} (x^2+\xi^2))
e^{-\frac{x^2+\xi^2}{h}} \ dx d\xi ,
$$
for a real constant $C$. Then applying the lemma with $F=1$ and $n=1$
proves that the integral term is zero.\\
If $j=k$ we obtain
$$
\int_{\R^2} W(\psi_j,\psi_j) (x,\xi) \ d\mu_{\R^2,h/2}(x,\xi)
=
\frac{(-1)^j}{\pi h} \int_{\R^2}  L_j^{(0)}( \frac{2}{h} (x^2+\xi^2))
e^{-\frac{x^2+\xi^2}{h}} \ dx d\xi .
$$
A polar change of variables and the expression of $ L_j^{(0)}$ then give
$$
\int_{\R^2} W(\psi_j,\psi_j) (x,\xi) \ d\mu_{\R^2,h/2}(x,\xi)
=
(-1)^j \int_{\R^+}  \sum_{m=0}^j \frac{j!}{(j-m)!(m!)^2}(-2u)^m e^{-u} \ du=1,
$$
using the fact that $\int_0^\infty u^m e^{-u} \ du = m!$.\\

As in the finite dimensional case, the interpretation is that
\begin{equation*}
\int_{\R^2} W(\psi_j,\psi_k) (x,\xi) \ d\mu_{\R^2,h/2}(x,\xi) =\langle
Op^{Weyl}_h(1) \psi_j,\psi_k
\rangle = \langle \psi_j,\psi_k
\rangle,
\end{equation*}
because $Op^{Weyl}_h(1)$ is the identity operator. 

\subsection{Decompositions over a Hilbert basis of $L^2(B,\mu_{B,h/2})$ }

Let $\hhh$ be a real, separable and infinite dimensional Hilbert space,
with an orthonormal basis   $(e_i)_{i\in \N^*}$.
Let $B$ be a Wiener extension of $\hhh$, for the measure $\mu_{B, h/2}$.\\

For a multiindex $\alpha$ we set
\begin{equation}\label{psialpha}
\forall y \in B,\quad 
\psi_{\alpha}^B(y) = \prod_{j\in \N^*} \psi_{\alpha_j}(\ell_{e_j}(y)).
\end{equation}
Although it runs over $\N^*$, this product is finite since the $\alpha_j$
are all $0$ but for a finite number and the factor $\psi_{\alpha_j}$  is equal to
$1$ if $\alpha_j=0$. The same formula defines Hermite functions of
finite dimension $d\geq 1$. When $\alpha$ runs over all multiindices, the family of the $\psi_{\alpha}^B$ is an orthonormal family of $L^2(B,\mu_{B,h/2})$
(see \cite{Jan}, Th. 2.6 and 3.21). \\

We get
\begin{lemm}\label{WignerB} 
For all multiindices  $\alpha$ and $\beta$, one has
$$
\forall (z,\zeta)\in B^2,\quad 
W_{h,B}(\psi_{\alpha}^B,\psi_{\beta}^B)(z,\zeta)
=
\prod_{j\in \N^*}W_{h,\R}(\psi_{\alpha_j},\psi_{\beta_j})
(\ell_{e_j}(z), \ell_{e_j}(\zeta)),
$$
with a mock infinite product once again. 
\end{lemm}
\noindent{\bf Proof } 

We may write 
$W_{h,B}(\psi_{\alpha}^B,\psi_{\beta}^B)(z,\zeta)$ as a finite dimensional
integral, using the distribution of the random vector $(\ell_{e_i})_{i}$,
the index running on all indices such that $\alpha_i$ or $\beta_i$ is not
equal to $0$.  This integral
splits into integrals over $\R^2$, since the distribution is normal
with diagonal covariance matrix and the integrated function itself is a product.
\hfill $\square$\\ 

When we use this lemma later on, we will express each  Wigner function  of
dimension $1$ in terms of the
Laguerre polynomials, according to \eqref{Wig-dim1}.\\

Recall that Formula \eqref{13-AJNJFA} defines a bilinear form
$Q_h^{Weyl}(\tilde{F})$ associated with a
function $\tilde{F}$,
defined on $B^2$. In certain cases,  it may be linked with an operator. In
view of a decomposition on the orthonormal  basis
$(\psi_{\alpha}^B)_{\alpha}$, we compute the 
$Q^{Weyl}_h(\tilde{F}) (\psi_{\alpha}^B,\psi_{\beta}^B)$.
For two multiindices $\alpha, \beta$, we consequently set
\begin{equation}\label{Ialphabeta-def}
I_{\alpha,\beta}( \tilde{F}) : =
\int_{B^2} \tilde{F}(z,\zeta) W_{h,B}(\psi_{\alpha}^B,\psi_{\beta}^B)(z,\zeta)
\  d\mu_{B^2, h/2}(z,\zeta) =Q^{Weyl}_h(\tilde{F}) (\psi_{\alpha}^B,\psi_{\beta}^B).
\end{equation}
 We address the case 
$\alpha\neq \beta$ just below. We  treat the case when $\alpha=\beta$,
under stronger conditions,
in Prop. \ref{prod-rad-psi}.

\begin{prop}\label{IPP-rad}
Suppose that $\tilde{F}$ is  cylindrical, based on 
$E={\rm Vect}(e_1,\dots, e_d)\subset B'$ and has the expression
\begin{equation}\label{expr-F}
\forall (z,\zeta)\in B^2, \quad \tilde{F}(z,\zeta)=
\Ddot{F}(\ell_{e_1}(z),\dots, \ell_{e_d}(z),\ell_{e_1}(\zeta),\cdots,
\ell_{e_d}(\zeta)),
\end{equation}
with $\Ddot{F}$ smooth ($C^{\infty}$) on $\R^{2d}$, increasing at most
polynomially, along with all its partial derivatives.\\
Consider $I_{\alpha,\beta}(\tilde{F})$ for  $\alpha \neq \beta$.\\
If  $\alpha_j\neq \beta_j$ for an index $j>d$,
then $I_{\alpha,\beta}(\tilde{F})=0$.\\
If, for all $j>d$,   $\alpha_j= \beta_j$, take $j\leq d$ such that
$\alpha_j\neq \beta_j$.
For any differentiation order $n$, one may write:
$$
I_{\alpha,\beta}(\tilde{F}) =
\frac{i ^n}{(\beta_j-\alpha_j)^n}
\int_{\R^{2d}} \left( (x_j\partial_{\xi_j}-\xi_j\partial_{x_j})^n \Ddot{F}\right)
(x,\xi) \prod_{l=1}^{d} W_{h,\R}(\psi_{\alpha_l},\psi_{\beta_l}) (x_l,\xi_l)
d\mu_{\R^d, h/2}(x,\xi).
$$
In particular, if
$ (x_j\partial_{\xi_j}-\xi_j\partial_{x_j})\Ddot{F}$
vanishes for all $j$,   $I_{\alpha,\beta}(\tilde{F})=0$ if   $\alpha \neq \beta$. 
\end{prop}

The proof of Prop. \ref{IPP-rad}
  is a consequence of the following lemma and of Lemma \ref{IPP}
in dimension $1$.

\begin{lemm}\label{Ialphabeta}  
Under the conditions and with  the notations of Prop.\ref{IPP-rad}  above,
if, for all $j>d$, $\alpha_j=\beta_j$, then the integral
$I_{\alpha,\beta}(\tilde{F})$ satisfies:
$$
I_{\alpha,\beta}(\tilde{F}) =\int_{\R^{2d}} \Ddot{F}(x,\xi) \ \prod_{j=1}^d
W_{h,\R}(\psi_{\alpha_j},\psi_{\beta_j})(x_j,\xi_j) \ 
e^{-\frac{1}{h}(|x|^2 + |\xi|^2)} \frac{dx d\xi}{(\pi h)^d}.
$$
If there exists $j>d$ such that
$\alpha_j\neq \beta_j$, then  $I_{\alpha,\beta}=0$.
\end{lemm}
\noindent{\bf Proof of the lemma}\\
Let $n$ be the largest index for which $\alpha_n$ or $\beta_n$
is not equal to $0$.\\
Suppose $n\leq d$. In this case, $\alpha_j=\beta_j=0$ for all $j> d$.
The measure on  $B$ (and on $B^2$)
decomposes as a product of Gaussian measures on  $E\times E^{\bot}$,
with $E^{\bot}= \{ x\in B, \ \forall u\in E, \ u(x)= 0 \}$,
as in \cite{RA}. This decomposition requires  that $E$ be a subset of $B'$
and not a more general subset of $\hhh$.
Every element $z$ of $B$ writes uniquely as  $z=z_E + z_{\bot}$, with 
$z_E=\sum_{i=1}^d e_i(z)e_i \in E$ and $u(z_{\bot}) =\ell_u(z_{\bot})= 0$ for all
$u\in E$.\\
Then, thanks to Lemma \ref{WignerB} and to Fubini's Theorem, we get,
by Formula \eqref{Ialphabeta-def}
$$
\begin{array}{lll}
I_{\alpha,\beta}(\tilde{F})
&=&\displaystyle 
\int_{E^2} \Ddot{F}(\ell_{e_1}(z_E),\dots, \ell_{e_d}(z_E),\ell_{e_1}(\zeta_E),\cdots,
\ell_{e_d}(\zeta_E))
\prod_{j=1}^dW_{h,\R}(\psi_{\alpha_j},\psi_{\beta_j})
(\ell_{e_j}(z_E), \ell_{e_j}(\zeta_E))\   d\mu_{E^2, h/2}\\ \\
& & \displaystyle \times 
\int_{(E^{\bot})^2}  1 
\  d\mu_{(E^{\bot})^2, h/2}.\\ \\ 
\end{array}
$$
The last integral is equal to $1$, the first one is equal to
$$
\int_{\R^{2d}} \Ddot{F}(x,\xi)
\prod_{k=1}^dW_{h,\R}(\psi_{\alpha_k},\psi_{\beta_k})
(x,\xi) \  d\mu_{\R^{2d}, h/2}(x,\xi),
$$
which gives the result. \\

Now suppose $n>d$. Let 
$G={\rm Vect}(e_{d+1},\dots, e_n)$.
 We decompose $B$  in the product of  $B=E\times G\times
(E\oplus G)^{\bot}$. 
The integral  $I_{\alpha,\beta} (\tilde{F})$ is  a product of three integral factors.
\begin{itemize}
\item
The last one is an integral on 
$( (E\oplus G)^{\bot})^2$ and it is equal to $1$ as in the preceding case.
\item
The first one is an integral on $E^2$ and has the same shape as in the
preceding case. 
\item
The integral in the middle is an integral on $G^2$, in which the integrated
function is only expressed thanks to Wigner functions. Indeed, $\tilde{F}$ and
$\Ddot F$ do not depend on the variables in $G$. This integral is then equal to 
a product of factors like 
$\int_{\R^2} W(\psi_{\alpha_s},\psi_{\beta_s}) (x,\xi) \ d\mu_{\R,h/2}$,
which are equal to
$\delta_{\alpha_s,\beta_s}$ by \eqref{Wigner-delta}.
\end{itemize}
If there exists $j>d$ such that $\alpha_j\neq \beta_j$, the product
$\prod_{d+1}^n \delta_{\alpha_s,\beta_s} $ is equal
to $0$ since this $j\leq n$.
If not, the integral ``in the middle'' is equal to $1$ and there
just remains the first factor. \\
This achieves the proof of Lemma \ref{Ialphabeta}\hfill $\square$
\\

\noindent{\bf Proof of Proposition  \ref{IPP-rad}}\\
When there is an index $j>d$ for which $\alpha_j\neq \beta_j$,
the proposition is a direct consequence of Lemma  \ref{Ialphabeta}.\\
Otherwise, take $j\leq d$ such that $\alpha_j\neq \beta_j$.
The Wigner function  corresponding to this index $j$ writes
$$
W_{h,\R}(\psi_{\alpha_j},\psi_{\beta_j})(x_j,\xi_j)=C_j
(x_j+ \eps i\xi_j)^{|\alpha_j-\beta_j|}
L^{|\alpha_j-\beta_j|}_{\min(\alpha_j,\beta_j)}(\frac{2}{h}(x_j^2+\xi_j^2)),
$$
with  $\eps =- {\rm sgn}(\alpha_j-\beta_j)$,  and the constant $C_j$ given by \eqref{Wig-dim1}. 
Set
$s=|\alpha_j-\beta_j|$. Applying Lemma \ref{IPP} with  a differentiation order
$n\in \N^*$ to the
integral on $(x_j,\xi_j)$ within the integral
$$
\int_{\R^{2d}} \Ddot{F}(x,\xi)
\prod_{k=1}^dW_{h,\R}(\psi_{\alpha_k},\psi_{\beta_k})
(x,\xi) \  d\mu_{\R^{2d}, h/2}(x,\xi)
$$
gives the result. 
This achieves the proof of Proposition \ref{IPP-rad}.\hfill $\square$

\subsection{Non positivity of the calculus }

This short paragraph contains the {\it proof of Proposition \ref{non-pos}},
which states a non positivity result analogous to the well-known result in
the finite dimensional case : an operator
with positive symbol is not necessarily positive. \\

The proof strongly relies on the finite dimensional situation. We choose
a symbol and a test function which give the result for the phase space $\R^2$.
We then build a cylindrical symbol and a cylindrical test
function adapted to our purpose. The computations are then exactly the same
as in the finite dimensional case,
for integrating  cylindrical functions gives rise to finite dimensional
integrals.\\

For  $a\in B'$ different from $0$ and $\nu >0$, one defines $F$ on $\hhh^2$
by
$$
\forall (x,\xi)\in \hhh^2,\quad F(x,\xi)= e^{-\nu ((a\cdot x )^2 + (a\cdot \xi)^2)}.
$$
Set $e_1= a/|a|$ and let  $(e_j)_{j\geq 1}$ be a Hilbert basis of $\hhh$,
consisting  of elements of $B'$ and  beginning with $e_1$.
Let $\Ddot{F}$ be the function defined on  $\R^2$ by
$\Ddot{F}(x,y)= e^{-\nu|a|^2(x^2+y^2)}$. It is rapidly decreasing,
hence $\gamma_{\R^2,s}\Ddot{F}$  is rapidly decreasing too for any
variance parameter $s>0$. According to Corollary  \ref{enfin}, the function $F$
admits, as a  stochastic extension, 
the function $\tilde{F}$ given by
$$
\tilde{F}(y,\eta) = e^{-\nu (\ell_a(y)^2+\ell_a(\eta)^2)} ,\ (y,\eta)\in B^2.
$$
This holds for any variance $s$.
\\
On the other hand, according to Lemma \ref{dans-SBeps},
$F\in S_m(\bbb,\eps)$ for every $m\in \N$, with respect to the sequence 
$\eps=\left(\frac{1}{n^2}\right)_{n\geq 1}$.
Taking $m=2$ one may, by Theorem \ref{1.4}, associate, with $F$, an operator 
$Op_h^{Weyl}(F)$ which is bounded on
$L^2(B,\mu_{B,h/2})$, for any $h>0$. Moreover, since $\ell_a\in\ddd_{B',h/2}$,
  \eqref{(1.17)} shows that the operator satisfies:
$$
\langle  Op_h^{Weyl}(F)\ell_a,\ell_a\rangle_{ L^2(B,\mu_{B,h/2})} = 
Q_h^{Weyl}(\tilde{F})(\ell_a,\ell_a)= \int_{B^2} \tilde{F}(y,\eta)
W_{h,B}(\ell_a,\ell_a)(y,\eta) \ d\mu_{B^2,h/2}(y,\eta).
$$
One has 
$\ell_a= |a| \sqrt{\frac{h}{2}} \psi_1(\ell_{e_1})$ by \eqref{Herm-prem},
then \eqref{Wig-prem}
gives $W_{h,\R}(\psi_1,\psi_1)$.
By Def. \ref{Wigner-infini}, we then have
$$W_{h,B}(\ell_a,\ell_a)(y,\eta) =|a|^2\frac{h}{2}
W_{h,\R}(\psi_1,\psi_1)(\ell_{e_1}(y),\ell_{e_1}(\eta))= 
|a|^2( - \frac{h}{2}  + (\ell_{e_1}(y)^2 + (\ell_{e_1}(\eta)^2 )).
$$

Consequently, we may write that 
$$
Q_h^{Weyl}(\tilde{F}))(\ell_a,\ell_a) =  |a|^2
\int_{B^2}  e^{-\nu |a|^2 (\ell_{e_1}(y)^2+\ell_{e_1}(\eta)^2)}
\left(\ell_{e_1}(y)^2 + \ell_{e_1}(\eta)^2 -\frac{h}{2}\right)
\ d\mu_{B^2, h/2} (y,\eta).
$$

The random vector
 $(\ell_{e_1}(y), \ell_{e_1}(\eta))$ is normally distributed with distribution
 $\nnn(0, \frac{h}{2} I_2)$. Therefore
$$
Q_h^{Weyl}(\tilde{F}))(\ell_a,\ell_a) =
\int_{\R^2}  e^{-\nu |a|^2 (u^2+v^2)}  |a|^2
\left(u^2 + v^2 -\frac{h}{2}\right)
\ e^{-(u^2+v^2)/h} (\pi h)^{-1} \ du dv
.$$
A polar decomposition then gives 
$$
\langle  Op_h^{Weyl}(F)\ell_a,\ell_a\rangle_{ L^2(B,\mu_{B,h/2})} = 
Q_h^{Weyl}(\tilde{F}))(\ell_a,\ell_a) =
\frac{h|a|^2}{2 (1+h \nu |a|^2)^2} (1-h \nu |a|^2)
.$$
Since this expression is negative for sufficiently large $\nu >0$, the
operator $ Op_h^{Weyl}(F)$ is not positive.
This concludes the proof of Proposition \ref{non-pos}.
 \hfill $\square$\\


\subsection{The Flandrin conjecture for  infinite dimensional Wigner functions}
\label{sec-Flandrin}

This conjecture, emitted   by Flandrin in 1988 in the article
\cite{Fl} (section $5$), concerns maximization and localization
of a signal in time and frequency.
The question, which remained open a long time, is to know whether it is
true  that, for any convex and bounded set $C$ of $\R^{2n}$
and any rapidly decreasing function $u$, the Wigner function
of $(u,u)$, normalized by
\begin{equation}\label{Wigner-DDL}
W_{}(u,v)(x,\xi) =\int_{\R} e^{-2i\pi z\xi} u(x+\frac{z}{2})
\bar{v}(x-\frac{z}{2}) \ dz,
\end{equation}
satisfies
$$
\int_{C}   W_{}(u,u)(x,\xi) dx d\xi \leq ||u||^2_{L^2(\R,dx)}.
$$
The result is true for dimension $2$ disks and Euclidean balls in more
general dimension. \\

In \cite{DDL} (Theorem 1.2 or p 31 of the article), the authors
consider the convex set $[0,a]^2$ or $(\R^+)^2$ (in this case,
we agree that $a=\infty$).
They prove that, for
$a>0$ sufficiently large or infinite, there exists, to the contrary,
a rapidly decreasing function $u_a$ such that
\begin{equation}\label{ineq-DDL}
\int_{[0,a[^2}   W_{}(u_a,u_a)(x,\xi) dx d\xi > ||u_a||^2_{L^2(\R,dx)}.
\end{equation}
The proof is complex and intricate. The
function $u_a$ is not explicit, for instance, and has no reason to
be ``simple''.\\

This beautiful result is easy to ``translate `` to the infinite
dimensional case. As in the preceding paragraph, we derive, from
the finite dimensional symbol and test function, a cylindrical symbol and
a cylindrical  test function.\\

We now turn to the {\it proof of Proposition \ref{Flandrin1}.}\\
Let $e_1\in B'$ such that $|e_1|_{\hhh}=1$.
Let $a$ be infinite or sufficiently large for  \eqref{ineq-DDL}
to hold.
Define  $\tilde{v_a}$ by
$$
\forall z \in B,\quad \tilde{v_a}(z)= v_a(\ell_{e_1}(z)),\quad {\rm with}\quad 
\gamma_{\R,h/2}v_a= u_a,
$$
which means explicitly that
$\forall x \in \R,\ v_a(x)= (\pi h)^{1/4} e^{x^2/2h} u_a(x)$.\\
Recall that
$\tilde{F_a}$ is given by
$$
\tilde{F_a}(z,\zeta) = \1_{[0,a[}(\ell_{e_1}(z))
\1_{[0,2\pi h a[}(\ell_{e_1}(\zeta)) .
$$

We have to prove that
$$
Q_h^{Weyl}(\tilde{F_a})(\tilde{v_a},
\tilde{v_a})=
 \int_{B^2} \tilde{F}_a(z,\zeta) W_{B,h}(\tilde{v_a},
\tilde{v_a})(z,\zeta)
\ d\mu_{B^2,h/2}(z,\zeta) >
||\tilde{v_a}||^2_{L^2(B, \mu_{B,h/2})}    .
$$
Remark that the use of the quadratic form  \eqref{13-AJNJFA} is licit because
$\tilde{F_a}$ is a bounded Borel function.\\
Since $u_a\in \sss(\R)$, $ \tilde{v_a}\in \ddd_{E,h/2}$ with
$E={\rm Vect}(e_1)$ and
$W_{h,B}(\tilde{v_a},
\tilde{v_a})(z,\zeta)= W_{h,E}({v_a},{v_a})( \ell_{e_1}(z), \ell_{e_1}(\zeta))$
by Def. \ref{Wigner-infini}.
This gives
$$
\begin{array}{lll}
Q_h^{Weyl}(\tilde{F_a})(\tilde{v_a},\tilde{v_a})&\displaystyle
= \int_{B^2} \1_{[0,a[}(\ell_{e_1}(z))\1_{[0,2\pi h a[}(\ell_{e_1}(z))
W_{h,\R}({v_a},{v_a})(\ell_{e_1}(z),\ell_{e_1}(\zeta)) \ d\mu_{B^2,h/2}(z,\zeta)\\ \\
&\displaystyle = \int_{\R^2} \1_{[0,a[}(x)\1_{[0,2\pi h a[}(\xi)
 W_{h,\R}({v_a},{v_a})( x,\xi) \
 e^{-\frac{1}{h} (x^2+\xi^2)}  \frac{1}{\pi h} \ dx d \xi
\end{array}
$$
because the random vector 
 $( \ell_{e_1}(z), \ell_{e_1}(\zeta))$ has the normal distribution
 $\nnn(0,\frac{h}{2} I_2)$.\\
Now, one can check that, for two functions $u$ and $v$
defined on $\R$ and such that 
$\gamma_{\R,h/2} u,\gamma_{\R,h/2} v $ are rapidly decreasing, one has
$$
e^{-\frac{1}{h} (x^2+\xi^2)} W_{h,\R}(u,v)(x,\xi)
=\frac{1}{2} W_{} (\gamma_{\R,h/2} u,\gamma_{\R,h/2} v)(x, \frac{\xi}{2\pi h}).
$$
Hence
$$
	 Q_h^{Weyl}(\tilde{F_a})(\tilde{v_a},\tilde{v_a})=  \int_{\R^2} \1_{[0,a[}(x)\1_{[0,2\pi h a[}(\xi)
 W_{} (\overbrace{\gamma_{\R,h/2} v_a}^{u_a},\gamma_{\R,h/2} v_a)(x, \frac{\xi}{2\pi h})
\  \frac{1}{2\pi h} dzd\zeta.
$$
It remains to set $\eta=  \frac{\xi}{2\pi h}$ and to exploit the results of
\cite{DDL} recalled above to obtain that 
$$
 Q_h^{Weyl}(\tilde{F_a})(\tilde{v_a},\tilde{v_a}) > ||u_a||_{L^2(\R,dx)} = ||v_a||^2_{L^2(\R, \mu_{\R,h/2})} =
||\tilde{v_a}||^2_{L^2(B, \mu_{B,h/2})}    .
$$
This achieves the proof of Proposition \ref{Flandrin1}.
\hfill $\square$\\


\subsection{Positivity for a symbol with radial properties}
\label{sec-posit-rad}
  
In this part, we aim at proving  Proposition  \ref{prod-rad-gene}
below, which generalizes  Proposition  \ref{posit-rad-gene} stated in
the introduction. This result concerns a tensor product
of cylindrical radial functions and not a strictly radial function. \\

Let us first state the hypotheses. 
\begin{description}
\item[H1]
The symbol $F$ is  {\it cylindrical.}\\
Let $\Ddot{F}$ be smooth on   $\R^{2d}$, bounded as well as all its partial
derivatives of any order. 
Let $(e_1,\dots, e_d)$ be an orthonormal family of $\hhh$ with vectors in $B'$
and let $E={\rm Vect}(e_1,\dots, e_d)$.\\
One defines $\tilde{F}$ on  $B^2$ and $F$ on  $\hhh^2$, as in \eqref{expr-F},
by
\begin{equation}\label{expr-F-bis}
\begin{array}{lll}
 \displaystyle
\forall (x,\xi)\in \hhh^2, \quad {F}(x,\xi)=
\Ddot{F}(e_1\cdot x, \dots, e_d \cdot x,e_1\cdot \xi , \dots, e_d \cdot \xi),\\
\displaystyle
\forall (z,\zeta)\in B^2, \quad \tilde{F}(z,\zeta)=
\Ddot{F}(\ell_{e_1}(z),\dots, \ell_{e_d}(z),\ell_{e_1}(\zeta),\cdots,
\ell_{e_d}(\zeta)).
\end{array}
\end{equation}
\item[H2]
The symbol $F$ is a {\it tensor product of radial functions}. \\
Precisely, let  $\{1,\dots, d\}$ split into $s$ pairwise disjoint parts
 $ D_j$ with cardinal $d_j >0$, $1\leq j\leq s$
(with $\sum_{j=1}^s d_j= d$). Set $d_0=0$ and suppose that
$$
\{1,\dots, d\} = \cup_{j=1}^s D_j, \quad {\rm with }\ \  D_j= \{d_0 + \dots + d_{j-1}+1,
\dots, d_0 + \dots + d_{j-1}+d_j\}.
$$
Denote by $x_{D_j}$ the variable corresponding to the coordinates
indexed by $D_j$, and adopt the same conventions for the dual variable $\xi$
and the couples $(x,\xi)$.
Suppose that

$$
\Ddot{F}(x,\xi) = \prod_{j=1}^s \Ddot{F}_j(x_{D_j},\xi_{D_j}), { \rm with} \ \
\Ddot{F}_j(x_{D_j},\xi_{D_j})= \Phi_j(|x_{D_j}|^2 +|\xi_{D_j}|^2),
$$
where the $\Phi_j$  are smooth on $\R^+$, satisfy 
$\Phi'_j \geq 0$ and are such that the 
$\Ddot{F}_j$ are smooth, bounded and with bounded derivatives of arbitrary
order. 
Here, $|\ |$ denotes the Euclidean norm on $\R^{d_j}$. 
\end{description}

Under these conditions,  the following result holds.
\begin{prop}\label{prod-rad-gene}
Suppose that $\tilde{F}$ and $\Ddot{F}$ satisfy conditions 
H1 and H2.\\
There exists an operator $Op^{Weyl}_h(F)$ associated with $F$ and bounded on 
  $L^2(B,\mu_{B,h/2})$.\\
Let $\tilde{f}\in L^2(B,\mu_{B,h/2})$. Then 
$$
\langle Op^{Weyl}_h(F) \tilde{f},\tilde{f}  \rangle_{L^2(B,\mu_{B,h/2})}    
\geq   \prod_{j=1}^s \left( \frac{1}{h} \int_0^{\infty}
\Phi_j(t)e^{-t/h} \ dt \right) \ ||\tilde{f}||^2_{L^2(B,\mu_{B,h/2})}.
$$
Finally, if $\tilde{f}$ is $\in \ddd_{B',h/2}$, this expression can be written
with the quadratic form:
$$
\langle Op^{Weyl}_h(F) \tilde{f},\tilde{f}  \rangle_{L^2(B,\mu_{B,h/2})}    =
\int_{B^2} \tilde{F}(z,\zeta) W_{h,B}(\tilde{f},\tilde{f})(z,\zeta)
\  d\mu_{B^2, h/2}(z,\zeta) .
$$
\end{prop}

An important argument in the  proof of this proposition
is the decomposition of  $\tilde{f}$ on the basis of the 
$\psi_{\alpha}^B$ defined in \eqref{psialpha}.
Then we are  brought back to integrals
like
$$
I_{\alpha,\beta}( \tilde{F})=
\int_{B^2} \tilde{F}(z,\zeta) W_{h,B}(\psi_{\alpha}^B,\psi_{\beta}^B)(z,\zeta)
\  d\mu_{B^2, h/2}(z,\zeta).
$$
The study of such integrals began in Proposition
\ref{IPP-rad}, which implies that, in the present case,
$I_{\alpha,\beta}( \tilde{F})=0$ if $\alpha\neq \beta$. It remains to study the
terms for which $\alpha = \beta$. \\
To that aim,  we  require a stronger condition, in order to apply the
result established in \cite{AJN-low}:
namely, the functions
$\Ddot{F}_j$ and their derivatives are supposed to be bounded (see H2).\\
In the article  \cite{AJN-low}, the fact that the symbol is radial
implies, as in  Proposition \ref{IPP-rad},  that the oblique terms (for
$\alpha \neq \beta$) vanish. But it is crucial too for the terms for which
$\alpha = \beta$, because it allows the use
of a radial change of variables. Nevertheless, under the present weaker
hypotheses of tensorization, it is possible to apply the result of
\cite{AJN-low} to the individual factors associated with the parts $D_j$.
Hence the result in  \cite{AJN-low} could be slightly generalized to
tensor products of radial functions.
\\

Proposition \ref{prod-rad-gene} relies on the following  result,  which
deals with a
couple of functions of the basis $(\psi_{\alpha}^B)_{\alpha}$, in view of
the said decomposition.
\begin{prop}  \label{prod-rad-psi}
Under the preceding hypotheses and with the notations above, for every
multiindex $\alpha$,
$$
\int_{B^2} \tilde{F}(z,\zeta) W_{h,B}(\psi_{\alpha}^B,\psi_{\alpha}^B)(z,\zeta)
\  d\mu_{B^2, h/2}(z,\zeta)
\geq   \prod_{j=1}^s \left( \frac{1}{h} \int_0^{\infty}
\Phi_j(t)e^{-t/h} \ dt \right),
$$
where $s$ is  the number of radial factors in the tensor product.
\\
In particular, if the  $\Phi_j$ are nonnegative,  this  expression is
nonnegative.  The nonnegativity of
the integral is sufficient. \\
If $\alpha\neq \beta$, then
$$
\int_{B^2} \tilde{F}(z,\zeta) W_{h,B}(\psi_{\alpha}^B,\psi_{\beta}^B)(z,\zeta)
\  d\mu_{B^2, h/2}(z,\zeta)
=0.
$$
\end{prop}

\noindent{\bf Proof of  Proposition \ref{prod-rad-psi} }\\
Using Lemma \ref{Ialphabeta} and the fact that the  Wigner function is itself a
tensor product, one has 
$$
\begin{array}{lll}
  \displaystyle
  I_{\alpha,\alpha}(\tilde{F}) =
  \prod_{j=1}^s \int_{\R^{2d_j}} \Ddot{F_j}(x_{D_j},\xi_{D_j}) \
\left(\prod_{l\in D_j}
 W_{h,\R}(\psi_{\alpha_l},\psi_{\alpha_l})(x_l,\xi_l)\right) \ 
e^{-\frac{1}{h}(|x_{D_j}|^2 + |\xi_{D_j}|^2)} \frac{dx_{D_j} d\xi_{D_j}}{(\pi h)^{d_j}}.
\end{array}
$$
We then apply Theorem 1.1 of   \cite{AJN-low} to every integral factor
of the product. The relationship between the Wigner function of the present
paper and  the one of  \cite{AJN-low} is given by
\begin{equation*}
W_{h,\R^d}(u,v)(x,\xi) = \frac{1}{2^d}
e^{\frac{1}{h}(|x|^2 + |\xi|^2)} H_h(\gamma_{\R^d,h/2}u,\gamma_{\R^d,h/2}v)(x,\xi).
\end{equation*}
Therefore
$$
\begin{array}{lll}
  \displaystyle
   \int_{\R^{2d_j}} \Ddot{F_j}(x_{D_j},\xi_{D_j}) \
\left(\prod_{l\in D_j}
 W_{h,\R}(\psi_{\alpha_l},\psi_{\alpha_l})(x_l,\xi_l)\right) \ 
 e^{-\frac{1}{h}(|x_{D_j}|^2 + |\xi_{D_j}|^2)} \frac{dx_{D_j} d\xi_{D_j}}{(\pi h)^{d_j}}\\ \\
 =\displaystyle
 \int_{\R^{2d_j}} \Ddot{F_j}(x_{D_j},\xi_{D_j}) \
 H_h(f_j,f_j )(x_{D_j},\xi_{D_j})
 \frac{dx_{D_j} d\xi_{D_j}}{(2 \pi h)^{d_j}},
\end{array}
$$
where the  function $f_j$ is given by 
$
f_j = \gamma_{\R^{d_j},h/2}\left( \prod_{l\in D_j}\psi_{\alpha_l}\right).
$

Theorem  1.1 of  \cite{AJN-low} then states that
$$
\int_{\R^{2d_j}} \Ddot{F_j}(x_{D_j},\xi_{D_j}) \
 H_h(f_j,f_j )(x_{D_j},\xi_{D_j})
 \frac{dx_{D_j} d\xi_{D_j}}{(2 \pi h)^d} \geq
 \frac{1}{h} \int_0^{\infty} \Phi_j(t)e^{-t/h} \ dt ||f_j||^2_{L^2(\R^{d_j}, dx)}.
$$
Note that the norm of $f_j$  is the norm in a $L^2$ space for
the Lebesgue measure.

Now,  $\gamma_{\R^d,h/2}$ is  an isometric isomorphism between 
 $L^2(\R^d,\mu_{\R^d,h/2})$ and
$L^2(\R^d,dy )$ (for the Lebesgue measure), hence $f_j$ has the same norm as 
$ \prod_{l\in D_j}\psi_{\alpha_l}$  in $L^2(\R^{d_j}, \mu_{\R^{d_j},h/2})$,
which is $1$.\\
This gives the first result of the proposition.\\
If the multiindices are different, one applies Proposition 
\ref{IPP-rad}, which proves that the expression is equal to $0$. \\
This concludes the proof of Proposition
\ref{prod-rad-psi}\hfill $ \square$\\ 

\noindent{\bf Proof  of Proposition \ref{prod-rad-gene} }\\
Recall that, for all  $(x,\xi)\in \hhh^2$,
$
F(x,\xi)= \Ddot{F}( e_1\cdot x,\dots, e_d \cdot x,e_1\cdot \xi,
\dots, e_d\cdot \xi).
$\\
Since the function $\Ddot{F}$ is smooth, bounded and with bounded derivatives,
the function $F$ above  is in $S_m(B,\eps)$ for arbitrary large $m\geq 2$,
with  $\eps=\left(\frac{1}{n^2}\right)_{n\geq 1}$ by Lemma  \ref{dans-SBeps}.
Its
stochastic extension  in
$L^p(B^2,\mu_{B,s})$ is  $\tilde{F}$, defined by \eqref{expr-F-bis},
for any $p\in [1,\infty[$ and any $s>0$.
\\
According to  Theorem \ref{1.4}, there exists
an operator  $Op^{Weyl}_h(F)$ which is bounded on $L^2(B,\mu_{B,h/2})$.
It satisfies, for all couple  $(\tilde{f},\tilde{g}))$ of $ \ddd_{B',h/2}$ 
$$
\langle Op^{Weyl}_h(F) \tilde{f},\tilde{g}\rangle_{L^2(B,\mu_{B,h/2})}
=
\int_{B^2} \tilde{F}(z,\zeta) W_{h,B}(\tilde{f},\tilde{g})(z,\zeta)
\ d \mu_{B,h/2}.
$$

Let $\tilde{f} \in L^2(B,\mu_{B,h/2})$ and decompose $\tilde{f}$
on the basis of the
$\psi_{\alpha}^B$. Set
$f_{\alpha}= \langle \tilde{f},\psi_{\alpha}^B \rangle e_{L^2(B,\mu_{B,h/2})}$.
Owing to the continuity of the operator, we obtain
$$
\langle Op^{Weyl}_h(F) \tilde{f},\tilde{f}\rangle_{L^2(B,\mu_{B,h/2})}
= \sum_{\alpha,\beta} \langle Op^{Weyl}_h(F) \psi_{\alpha}^B,
\psi_{\beta}^B\rangle_{L^2(B,\mu_{B,h/2})}
f_{\alpha}
\overline{f_{\beta}}.
$$

The functions $\psi_{\alpha}^B$ are in $\ddd_{B',h/2}$, hence the factors
$\langle Op^{Weyl}_h(F) \psi_{\alpha}^B,\psi_{\beta}^B\rangle_{L^2(B,\mu_{B,h/2})}$
can be expressed with the Wigner function. For $\alpha\neq \beta$
we saw that the factor is $0$
(Prop. \ref{IPP-rad} and \ref{prod-rad-gene}). It remains 
$$
\langle Op^{Weyl}_h(F) \tilde{f},\tilde{f} \rangle_{L^2(B,\mu_{B,h/2})}
 =\sum_{\alpha} f_{\alpha} 
\overline{f_{\alpha}}\ \times \ 
\int_{B^2} \tilde{F}(z,\zeta) W_{h,B}(\psi_{\alpha}^B,\psi_{\alpha}^B)(z,\zeta)
\  d\mu_{B^2, h/2}(z,\zeta).
$$
Now, all the integrals for  $\alpha=\beta$ are greater than 
$  \displaystyle \prod_{j=1}^s \left( \frac{1}{h} \int_0^{\infty}
\Phi_j(t)e^{-t/h} \ dt \right) $ , by Prop. \ref{prod-rad-psi}.
We then get
$$
\begin{array}{lll}
  \displaystyle 
\langle Op^{Weyl}_h(F) \tilde{f},\tilde{f}\rangle_{L^2(B,\mu_{B,h/2})}
&  \displaystyle \geq   \prod_{j=1}^s \left( \frac{1}{h} \int_0^{\infty}
\Phi_j(t)e^{-t/h} \ dt \right)  \sum_{\alpha} |f_{\alpha}|^2 \\ \\ &
 \displaystyle \geq \prod_{j=1}^s \left( \frac{1}{h} \int_0^{\infty}
\Phi_j(t)e^{-t/h} \ dt \right)\  ||\tilde{f}||_{{L^2(B,\mu_{B,h/2})}}^2.
\end{array}
$$
This concludes the proof of Proposition \ref{prod-rad-gene}.
\hfill $ \square$\\ 


\section{G\aa rding's inequality}\label{sec-Gaarding}

We must recall the definitions and results which have not taken place
in the introduction. In particular, we define the partial heat operators
which act on functions defined on the Wiener space as in  \cite{HA,RA}
and the hybrid quadratic forms associated with a symbol $F$ in
$S_m(\bbb,\eps)$ and a subspace $E$ of $\fff(B')$. These notions
allow us to give the decomposition  \eqref{dec-TS} of the symbol $F$.
The decomposition was a step in the construction of the operator
$Op^{Weyl}_h(F) $   
 in \cite{AJN-JFA}. Here it is useful to isolate a term which is nonnegative. 

 \subsection{Heat operators, hybrid and Anti-Wick quadratic forms}

The following notions are defined in 
 \cite{AJN-JFA} (sections 2 and 3).

Let $E\in \fff(B')$  and denote by  
$E^{\bot}$ the following subspace of $B$:
$$
E^{\bot}= \{ x\in B, \ \forall u\in E, \ u(x)= 0 \}.
$$
According to  \cite{RA}, $(E^{\bot}\cap \hhh, E^{\bot})$ is  a Wiener space and,
for all $h>0$,  $ \mu_{B,h} =\mu_{E,h}\otimes \mu_{E^{\bot},h}$.
The same decomposition is valid for $B^2$ and already appears in the proof of
Lemma \ref{Ialphabeta}.
We denote by  $ X_E, X_{E^{\bot}}$ the variables of $E^2$
and $(E^{\bot})^2$ which appear below, in the integrals. 
\\
For  $E\in \fff(B')$ one defines, for $t>0$, the following partial heat
operators. For every Borel, bounded function $\tilde{F}$ defined on $B^2$, we set
\begin{equation}\label{chal-part}
\begin{array}{lll}
\displaystyle
\forall X\in B^2, &\displaystyle
\widetilde{H_{E,t}}(\tilde{F})(X)  & =
\displaystyle
\int_{E^2} \tilde{F}(X+Y_E)\ d\mu_{E^2,t}(Y_E)\\
\displaystyle
\forall X\in B^2,  &\displaystyle
\widetilde{H_{E^{\bot}},t}(\tilde{F})(X)  & =
\displaystyle
\int_{{E^{\bot}}^2} \tilde{F}(X+Y_{E^{\bot}})\ d\mu_{{E^{\bot}}^2,t}(Y_{E^{\bot}}).\\ 
\end{array}
\end{equation}
The functions obtained are Borel functions, bounded on $B^2$. 
\\
Let $E_1\subset E_2\in \fff(B')$ and let $S\subset E_2$ be the orthogonal
complement of   $E_1$  in $E_2$.
The heat operators satisfy
\begin{equation}\label{E1E2S}
\widetilde{H_{E_1^{\bot},t}} = \widetilde{H_{E_2^{\bot},t}}\widetilde{H_{S,t}}.
\end{equation}
For a given, bounded, Borel function  $\tilde{F}$  defined on $B^2$ and a
subset $E\in \fff(B')$, let us define, as in 
\cite{AJN-JFA} (Definition 2.2), a hybrid quadratic form 
$Q_h^{hyb,E}(\tilde{F})$ on  $\ddd_{B',h/2}^2$, setting, for $\tilde{f},\tilde{g}$ in
$\ddd_{B',h/2}$ :
\begin{equation}\label{Qhyb}
Q_h^{hyb,E}(\tilde{F})(\tilde{f},\tilde{g})= Q_h^{Weyl}(\widetilde{H_{E^{\bot},h/2}}\tilde{F})(\tilde{f},\tilde{g}).
\end{equation}
Recall that $ Q_h^{Weyl}$ is defined by \eqref{13-AJNJFA}.

If $E_1\subset E_2\in \fff(B')$ and if $S\subset E_2$ is the orthogonal
complement of $E_1$ in  $E_2$, the hybrid forms satisfy the relationship
\begin{equation}\label{29}
Q_h^{hyb,E_1}(\tilde{F})  = Q_h^{hyb,E_2}(\widetilde{H_{S,h/2}}\tilde{F}).
\end{equation}

In the particular case when
$E=\{ 0\}$, $E^{\bot}=B$ and the heat operator is the classical heat operator
(classical in the theory of Wiener spaces : \cite{HA}). It allows one to
define a
quadratic form
$Q_h^{AW}(\tilde{F})$
on $\ddd_{B',h/2}^2$. This is
done in  Formula (23) in
\cite{AJN-JFA}:
\begin{equation}\label{23}
  \forall (\tilde{f},\tilde{g})\in \ddd_{B',h/2}^2,\quad
Q_h^{AW}(\tilde{F})(\tilde{f},\tilde{g}) = Q_h^{Weyl}(\widetilde{H_{B,h/2}}\tilde{F})(\tilde{f},\tilde{g}) .
\end{equation}
 In the finite dimensional case, the quadratic form
defined this way is linked
with the Anti-Wick operator with symbol $F$, hence the denomination
of the form introduced above.
According to Corollary 4.12 of \cite{AJN-JFA}, if $\tilde{F}$ is  
nonnegative, so is $Q_h^{AW}(\tilde{F})$, exactly as for the finite dimensional
Anti-Wick calculus:
$$
\forall \tilde{f} \in \ddd_{B',h/2},
\ Q_h^{AW}(\tilde{F})(\tilde{f},\tilde{f}) \geq 0.
$$

It remains to introduce operators which  give a decomposition of 
identity. For $j\geq 1$ we set 
$D_j={\rm Vect}( (e_j,0),(0,e_j)) \in \fff(B'^2)$.
Then, for any finite subset  $J$ of $\N^*$, we set:

\begin{equation}\label{T-S}
\widetilde{T_{J,h}}= \prod_{j\in J}(I-\widetilde{H_{D_j,h/2}}),\quad
\widetilde{S_{J,h}}= \prod_{j\in J}\widetilde{H_{D_j,h/2}},
\end{equation}
agreeing that, for $J=\emptyset$, the operator is $I$.

\subsection{Proof of G\aa rding's inequality}
This proof follows, on the one hand the arguments of 
\cite{AJN-JFA} which, in this infinite dimensional frame, 
define  $Op^{Weyl}_h(F)$ and, on the other hand, the method used in 
\cite{L-N} to establish a G\aa rding's inequality, in which the constant
does not depend on the dimension.  
We  highlight a term, corresponding to the Anti-Wick quadratic form,
which thus has positivity properties. \\

Let $\Lambda$ be a finite subset of $\N^*$. Set $E(\Lambda)=
{\rm Vect}(e_j,j\in \Lambda).$ One checks that
\begin{equation}\label{dec-TS}
\tilde{F} = \sum_{J\subset \Lambda} \widetilde{T_{J,h}}
\widetilde{S_{\Lambda\setminus J,h}}\tilde{F}.
\end{equation}
Hence
$$
Q^{hyb, E(\Lambda)}_h(\tilde{F}) =
\sum_{J\subset \Lambda}Q^{hyb, E(\Lambda)}_h(\widetilde{T_{J,h}}
\widetilde{S_{\Lambda \setminus J,h}}\tilde{F}).
$$
One isolates the term corresponding to the empty set  $J=\emptyset$, 
for which
$\widetilde{T_{\emptyset,h}}= I$.
This term is
$$
\begin{array}{llllll}
\displaystyle Q^{hyb, E(\Lambda)}_h(\widetilde{S_{\Lambda ,h}}\tilde{F})
& =  \displaystyle  Q^{hyb, E(\Lambda)}_h(\widetilde{H_{E(\Lambda) ,h/2}}\tilde{F})
\\ \\
& =  \displaystyle  Q^{hyb, \{0\}}_h(\tilde{F})  & {\rm by } \  \eqref{29}\
{\rm with }\ E_2=S=E(\Lambda)
\\ \\
& =  \displaystyle  Q^{Weyl}_h(\widetilde{H_{B,h/2}}\tilde{F})  & {\rm by } \  \eqref{Qhyb}\
\\ \\
& =  \displaystyle  Q^{AW}_h(\tilde{F})  & {\rm by } \  \eqref{23}.\
\\ \\
\end{array}
$$
Let $\tilde{f} \in \ddd_{B',h/2}$. Since $ Q^{AW}_h(\tilde{F})
(\tilde{f} ,\tilde{f} )\geq 0$,
$$
Q^{hyb, E(\Lambda)}_h(\tilde{F})(\tilde{f},\tilde{f}) \geq 
\sum_{J\subset \Lambda, J\neq \emptyset} Q^{hyb, E(\Lambda)}_h(\widetilde{T_{J,h}}
\widetilde{S_{\Lambda \setminus J,h}}\tilde{F})(\tilde{f},\tilde{f}).
$$
Now,
$$
Q^{hyb, E(\Lambda)}_h(\widetilde{T_{J,h}}\widetilde{S_{\Lambda \setminus J,h}}\tilde{F})(\tilde{f},\tilde{f})
=
Q^{hyb, E(J)}_h(\widetilde{T_{J,h}}\tilde{F})(\tilde{f},\tilde{f})
$$
according to \eqref{29} with  $E_2=E(\Lambda), S=E(\Lambda \setminus J)$ and $E_1=E(J)$.

Therefore, we obtain
$$
\begin{array}{llll}
\displaystyle
Q^{hyb, E(\Lambda)}_h(\tilde{F})(\tilde{f},\tilde{f})  &\displaystyle  \geq 
\sum_{J\subset \Lambda, J\neq \emptyset}
Q^{hyb, E(J)}_h(\widetilde{T_{J,h}}\tilde{F})(\tilde{f},\tilde{f})  \\ \\
  &\displaystyle  \geq -
\sum_{J\subset \Lambda, J\neq \emptyset}
\vert Q^{hyb, E(J)}_h(\widetilde{T_{J,h}}\tilde{F})(\tilde{f},\tilde{f})\vert  \\ \\
 &\displaystyle  \geq - M 
\sum_{J\subset \Lambda, J\neq \emptyset}
(81\pi h S_{\eps})^{|J|}\prod_{j\in J}\eps^2_j \ ||\tilde{f} ||^2_{L^2(B,\mu_{B,h/2})}
\end{array}
$$
using Formula  (35)  in Proposition 3.1
of \cite {AJN-JFA}, where $M= ||F||_{S_{2}(\bbb,\eps)}$ and $S_{\eps}$ is defined in
Prop. \ref{Gaa}.
Set $\lambda_j= 81\pi h S_{\eps}\eps^2_j $. The sequence
$(\lambda_n)_{n\geq 1}$ is summable, since $\eps$ is square summable.
One  has
$$
\sum_{J\subset \Lambda, J\neq \emptyset}\prod_{j\in J}\lambda_j
\leq \sum_{j\in \Lambda}\lambda_j \ \prod_{s\in \Lambda} (1+\lambda_s),
$$
(this is formula (20) of \cite{L-N}).
Hence, for every finite subset $\Lambda$ of $\N^*$, under the conditions
of Prop. \ref{Gaa}, one has, for
$\tilde{f} \in \ddd_{B',h/2}$
\begin{equation}\label{INEG}
Q^{hyb, E(\Lambda)}_h(\tilde{F})(\tilde{f},\tilde{f})\geq 
- ||F||_{S_{2}(\bbb,\eps)}  \sum_{j\in \Lambda}\lambda_j \ \prod_{s\in \Lambda} (1+\lambda_s)\ 
 ||\tilde{f} ||^2_{L^2(B,\mu_{B,h/2})}  .
\end{equation}

Let now $(\Lambda_n)_{n\in \N^*}$ be an increasing sequence of finite subsets of $\N^*$, with union equal to  $\N^*$. With every $n$, Proposition  3.2 of
\cite{AJN-JFA} associates an operator $Op_h^{hyb,E(\Lambda_n) }(\tilde{F})$, bounded in  $L^2(B,\mu_{B,h/2})$ and  satisfying 
$$
\forall (\tilde{f},\tilde{g})\in \ddd_{B',h/2}^2,\quad
Q^{hyb, E(\Lambda_n)}_h(\tilde{F})(\tilde{f},\tilde{g}) =
\langle Op_h^{hyb,E(\Lambda_n) }(\tilde{F})\tilde{f},\tilde{g} \rangle_{L^2(B,\mu_{B,h/2})}.
$$
Moreover, the sequence of these operators is a Cauchy sequence in 
the space $\lll(L^2(B,\mu_{B,h/2}))$  and converges to   the operator
$Op_h^{Weyl}(F)$. Originally, it is what defines $Op_h^{Weyl}(F)$ in
\cite{AJN-JFA}.
 Using \eqref{INEG} with  $\Lambda=\Lambda_n$, we get
$$
Op_h^{hyb,E(\Lambda_n) }(\tilde{F})(\tilde{f},\tilde{f})\geq 
- ||F||_{S_{2}(\bbb,\eps)}  \sum_{j\in \Lambda_n}\lambda_j \ \prod_{s\in \Lambda_n} (1+\lambda_s)\ 
||\tilde{f} ||^2_{L^2(B,\mu_{B,h/2})}  .
$$
Letting $n$ go to infinity yields
$$ Op_h^{Weyl}({F})(\tilde{f},\tilde{f})\geq 
- ||F||_{S_{2}(\bbb,\eps)}  \sum_{j\geq 0}\lambda_j \ \prod_{s\geq 0} (1+\lambda_s)\ 
||\tilde{f} ||^2_{L^2(B,\mu_{B,h/2})}  ,
$$
since the sequence  $(\lambda_n)_{n\geq 1}$ is  summable.\\
By density of $\ddd_{B',h/2}$  in $L^2(B,\mu_{B,h/2})$ and because of the continuity of $ Op_h^{Weyl}({F})$, this inequality still holds 
for an arbitrary
$\tilde{f} \in L^2(B,\mu_{B,h/2})$.
This concludes the proof of Proposition
 \ref{Gaa}.\hfill $\square$\\
 


\medskip

\noindent lisette.jager@univ-reims.fr\\
Universit\'e de Reims Champagne-Ardenne\\
Laboratoire de Math\'ematiques de Reims UMR CNRS 9008\\
Moulin de la Housse, BP 1039, 51687 Reims Cedex 2, France

\end{document}